\documentclass[12pt,letter,leqno]{article}
\usepackage[numbers]{natbib}
\usepackage{fancyheadings}
\usepackage[]{hyperref}
\usepackage{xcolor}
\hypersetup{colorlinks = true, citecolor = blue, linkcolor=blue, linkbordercolor = {white}}
\usepackage{amsmath}
\usepackage{amsfonts}
\usepackage{mathrsfs}
\usepackage{amsthm}
\usepackage{amssymb}
\usepackage{mathrsfs}
\usepackage{mathtools}
\usepackage{bm}
\usepackage[normalem]{ulem}
\usepackage{graphicx}
\usepackage{caption}

\newtheorem{theorem}{Theorem}
\newtheorem{lemma}{Lemma}
\newtheorem{corollary}{Corollary}
\newtheorem{proposition}{Proposition}

\newtheorem{example}{Example}
\newtheorem{remark}{Remark}

\newtheorem{innercustomgeneric}{\customgenericname}
\providecommand{\customgenericname}{}
\newcommand{\newcustomtheorem}[2]{%
  \newenvironment{#1}[1]
  {%
   \renewcommand\customgenericname{#2}%
   \renewcommand\theinnercustomgeneric{##1}%
   \innercustomgeneric
  }
  {\endinnercustomgeneric}
}
\newcustomtheorem{customlemma}{Lemma}
\newcustomtheorem{customproposition}{Proposition}
\newcustomtheorem{customtheorem}{Theorem}
\newcustomtheorem{customdefinition}{Definition}
\newcustomtheorem{customexample}{Example}

\newcommand{\lng}{{\langle}}
\newcommand{\rng}{{\rangle}}

\newcommand{\mbS}{\mathbb{S}}

\newcommand{\Bl}{ B^l}
\newcommand{\Blq}{ B^l_q}
\newcommand{\Blr}{ B^l_r}
\newcommand{\Bls}{ B^l_s}
\newcommand{\suml}{\sum_{l=0}^{[T_n]}}
\newcommand{\sumq}{\sum_{q=1}^{N(d,l)}}
\newcommand{\sumr}{\sum_{r=1}^{N(d,l)}}
\newcommand{\sums}{\sum_{s=1}^{N(d,l)}}

\newcommand{\malpha}{\bm{\alpha}}

\newcommand{\tran}{^{\mathstrut\scriptscriptstyle{\top}}}

\newcommand{\RE}{{\rm Re}}
\newcommand{\IM}{{\rm Im}}

\newcommand{\op}{{\rm op}}
\newcommand{\E}{{\rm E}}
\newcommand{\EL}{{\rm EL}}
\newcommand{\Var}{{\rm Var}}
\newcommand{\Leb}{{\rm Leb}}

\newcommand{\Tr}{{\rm Trace}}

\renewcommand{\theequation}{\thesection.\arabic{equation}}

\begin{document}

\Large
\centerline{\bf Density estimation and regression analysis on $\mathbb{S}^d$}
\centerline{\bf in the presence of measurement error}
\large
\vspace{.2cm} \centerline{Jeong Min Jeon and Ingrid Van Keilegom}
\vspace{.2cm} \centerline{Research Centre for Operations Research and Statistics}
\vspace{.2cm} \centerline{KU Leuven, Belgium}
\vspace{.2cm}
\normalsize
\begin{abstract}
This paper studies density estimation and regression analysis with contaminated data observed on the unit hypersphere $\mbS^d$ for $d\in\mathbb{N}$. Our methodology and theory are based on harmonic analysis on general $\mbS^d$. We establish novel nonparametric density and regression estimators, and study their asymptotic properties including the rates of convergence and asymptotic distributions. We also provide asymptotic confidence intervals based on the asymptotic distributions of the estimators and on the empirical likelihood technique. We present practical details on implementation as well as the results of numerical studies.
\end{abstract}
\small
\begin{quotation}
\end{quotation}
\normalsize

\section{Introduction}

\setcounter{equation}{0}
\setcounter{subsection}{0}

Statistical analysis with data involving measurement errors has been a challenging problem in statistics. When some variables are not precisely observed due to measurement errors, direct application of existing methods designed for error-free variables results in incorrect inference. To explain this, let us consider a simple case
where both the covariate $X$ and the response $Y$ are real-valued. To estimate the regression function $m(x)=\E(Y|X=x)$ at a point $x$, one may apply
`local smoothing' to $Y_i$ around each point $x$. For example, the Nadaraya-Watson estimator of $m$ is to take a weighted average of $Y_i$
corresponding to $X_i$ that fall in a neighborhood of each point $x$.
This makes sense since $Y_i$ corresponding to $X_i$ near $x$ have `correct' information about $m(x)$.
Now, suppose that $X_i$ are not available but $Z_i=X_i+U_i$ are where $U_i$ are unobserved measurement errors.
In this case, the naive approach, simply taking a weighted average of $Y_i$ corresponding to $Z_i$
that fall in a neighborhood of the point $x$, should fail since $X_i$ corresponding to such $Z_i$ may locate far away from $x$ and thus
the corresponding $Y_i$ may not have correct information about $m$ at $x$. To treat this issue, appropriate correction methods have been proposed. To list only a few, \cite{Stefanski and Carroll (1990)} introduced a deconvolution kernel density estimator,
and \cite{Fan (1991a)} and \cite{Fan (1991b)} studied its rate of convergence and asymptotic distribution, respectively.
Based on the deconvolution kernel, \cite{Fan and Truong (1993)} investigated the rate of convergence
for a Nadaraya-Watson-type regression estimator and \cite{Delaigle et al. (2009)} studied the asymptotic distribution
of a local-polynomial-type regression estimator. For an introduction to the measurement error problems, we refer to \cite{Meister (2009)} and \cite{Delaigle (2014)}. However, the aforementioned works are restricted to Euclidean data.

Analyzing non-Euclidean data is becoming an important topic in modern statistics due to rapidly emerging non-Euclidean data in various fields. It is challenging since there is no vector space structure on non-Euclidean spaces in general.
For a recent review on non-Euclidean data analysis, we refer to \cite{Marron and Alonso (2014)}. Data observed on the unit hypersphere $\mbS^d=\{x\in\mathbb{R}^{d+1}:\|x\|=1\}$ for $d\in\mathbb{N}$, called hyperspherical data, are one of the most abundant non-Euclidean data. Hyperspherical data include circular data ($d=1$), spherical data ($d=2$) and other higher dimensional data (e.g. \cite{Scealy and Welsh (2011)}, \cite{Garcia-Portugues et al. (2016)}). Previous works on error-free circular, spherical or general hyperspherical data include density estimation (\cite{Hall et al. (1987)}, \cite{Garcia-Portugues et al. (2013)}), regression analysis (\cite{Chang (1989)}, \cite{Rivest (1989)}, \cite{Rosenthal et al. (2014)}, \cite{Jeon et al. (2021)}) and statistical testing (\cite{Cuesta-Albertos et al. (2009)}, \cite{Boente et al. (2014)}, \cite{Garcia-Portugues et al. (2020)}). Among them, \cite{Chang (1989)} did not cover a measurement error problem and simply considered the case where both response and predictor are spherical variables and the response is symmetrically distributed around the product of an unknown orthogonal matrix and the predictor. For a recent review on hyperspherical data analysis, we refer to \cite{Pewsey and Garcia-Portugues (2021)}.

Some areas that hyperspherical data arise are meteorology and astronomy. However, data from such areas are prone to contain measurement errors due to the technical limitations of measuring devices. For example, measuring the exact wind direction, positions of sunspots on the sun or direction from the earth to an astronomical object is not easy since they move very fast and/or they are very far (e.g. \cite{Baranyi et al. (2001)}, \cite{Gao et al. (2015)}). Also, such measurements are sometimes disturbed by some substances in between. In addition, each observation vector in Euclidean data is sometimes normalized to have the unit norm to ensure that data analysis is only affected by the relative magnitudes of vector elements rather than the absolute magnitudes of vectors themselves. If the original Euclidean data contain measurement errors, then the resulting hyperspherical data also contain measurement errors.

In spite of the importance of analyzing contaminated hyperspherical data, there exist only few works, and most of the existing works are restricted to deconvolution density estimation on either $\mbS^1$ or $\mbS^2$ (e.g. \cite{Efromovich (1997)}, \cite{Healy et al. (1998)}, \cite{Kim and Koo (2002)}, \cite{Kim et al. (2004)}). Some other works for other types of contaminated non-Euclidean data include deconvolution density estimation on special orthogonal groups (\cite{Kim (1998)}), compact and connected Lie groups (\cite{Kim and Richards (2001)}), the Poincar\'{e} upper half plane (\cite{Huckemann et al. (2010)}) and the 6-dimensional Euclidean motion group (\cite{Luo et al. (2011)}). All the aforementioned works on deconvolution density estimation studied only the rates of convergence of their estimators. Recently, \cite{Jeon et al. (2022)} studied density estimation and regression analysis with contaminated Lie-group-valued predictor. To the best of our knowledge, \cite{Jeon et al. (2022)} is the unique work that considered regression analysis with contaminated manifold-valued variables. However, since $\mbS^d$ for $d=2$ and $d\geq4$ are not a Lie group, it is important to study such unexplored cases.

In this paper, our primary aim is to develop the a deconvolution regression estimator on $\mbS^2$ and investigate its rates of convergence. We also aim to construct the asymptotic distributions and asymptotic confidence intervals for both deconvolution density and regression estimators on $\mbS^2$. Those have not been studied in the literature  despite their importance. To achieve them in a more general setting, we instead study deconvolution density estimation and regression analysis on $\mbS^d$ for $d\in\mathbb{N}$. These general problems also have not been considered in the literature. Our deconvolution density estimator on $\mbS^d$ generalizes the deconvolution density estimator on $\mbS^1$ introduced in \cite{Efromovich (1997)} and the one on $\mbS^2$ introduced in \cite{Healy et al. (1998)}. Our deconvolution regression estimator on $\mbS^d$ also generalizes the deconvolution regression estimator on $\mbS^1$ introduced in \cite{Jeon et al. (2022)}. We build up a theoretical foundation for those general estimators. We establish several finite-sample properties of the estimators. We also study the uniform consistency of the density estimator and the rates of convergence for both density and regression estimators. In addition, we derive the asymptotic distributions and two types of asymptotic confidence intervals for both estimators under a high-level condition. The high-level condition is verified for certain cases. Moreover, we present several numerical studies and some practical details on implementation which have received less attention in the literature in spite of their importance. We emphasize that deriving the results in this paper is quite different from the ways in the Euclidean case since it is based on hyperspherical harmonic analysis which is less considered in statistics. Also, dealing with $\mbS^d$ is more challenging than dealing with $\mbS^2$ since general hyperspherical harmonic analysis is much more complex than harmonic analysis on $\mbS^2$. Indeed, it leads to more complex analysis for every result and requires broader discussions.

This paper is organized as follows. In Section \ref{methodology}, we introduce general hyperspherical harmonic analysis with some practical examples and our estimators with some finite-sample properties. The rates of convergence and asymptotic distributions of our estimators are shown in Section \ref{rate and distribution}. We construct the asymptotic confidence intervals in Section \ref{confidence interval}, and present the simulation studies and real data analysis in Section \ref{simulation}. The Supplementary Material contains additional practical details and all technical proofs. %We collect remaining technical proofs in the Supplementary Material.

\section{Preliminaries and methodology}\label{methodology}

\setcounter{equation}{0}
\setcounter{subsection}{0}

\subsection{Preliminaries}

Our methodology is largely based on harmonic analysis on the $d$-dimensional unit hypersphere $\mbS^d$ for $d\in\mathbb{N}$. Here, we give a brief introduction on it. Further details can be found in \cite{Atkinson and Han (2012)} and \cite{Efthimiou and Frye (2014)}.

A function $f:\mathbb{R}^{d+1}\rightarrow\mathbb{C}$ is called a harmonic homogeneous polynomial of degree $l\in\mathbb{N}_0:=\{0\}\cup\mathbb{N}$ in $d+1$ variables if $f$ takes the form
\begin{align*}
f(t_1,\ldots,t_{d+1})=\sum_{\malpha=(\alpha_1,\ldots,\alpha_{d+1})\in\mathbb{N}^{d+1}_0:\,\sum_{i=1}^{d+1}\alpha_i=l}c_{\malpha}\cdot\prod_{i=1}^{d+1}t_i^{\alpha_i}
\end{align*}
for $c_{\malpha}\in\mathbb{C}$ and satisfies $\sum_{i=1}^{d+1}\partial^2f(t_1,\ldots,t_{d+1})/\partial t_i^2\equiv0$. For such $f$, the domain restricted function $f|_{\mbS^d}:\mbS^d\rightarrow\mathbb{C}$ is called a spherical harmonic of order $l$ in $d+1$ variables. We denote the space of all spherical harmonics of degree $l$ in $d+1$ variables by $\mathfrak{B}^l(\mbS^d)$ and call it the spherical harmonic space of order $l$ in $d+1$ variables.

It is known that $\mathfrak{B}^l(\mbS^d)$ is a vector space of dimension
\begin{align*}
N(d,l)=\frac{(2l+d-1)\cdot(l+d-2)!}{l!\cdot(d-1)!}.
\end{align*}
Direct computations show that $N(1,l)=2$ and $N(2,l)=2l+1$ for $l\in\mathbb{N}$, and $N(d,0)=1$. It is also known that the vector space spanned by $\{\mathfrak{B}^l(\mbS^d):l\in\mathbb{N}_0\}$ is a dense subspace of the $L^2$ space $L^2((\mbS^d,\nu),\mathbb{C})$, where $\nu$ is the scaled spherical measure on $\mbS^d$ defined by \begin{align*}
\nu(A)=\frac{{\rm{area}}(\mbS^d)}{\Leb(B(0,1))}\cdot\Leb(\{ta:t\in[0,1],a\in A\})
\end{align*}
for any Borel subset $A$ of $\mbS^d$, where ${\rm{area}}(\mbS^d)$ is the surface area of $\mbS^d$, $\Leb$ is the Lebesgue measure on $\mathbb{R}^{d+1}$ and $B(0,1)$ is the closed ball centered at zero with radius one. We note that $\nu(\mbS^d)={\rm{area}}(\mbS^d)$.
We also note that $L^2((\mbS^d,\nu),\mathbb{C})$ is a separable Hilbert space with inner product $\lng f,g\rng_2=\int_{\mbS^d}f(x)\overline{g(x)}\,d\nu(x)$, where $\overline{g(x)}$ is the conjugate of $g(x)$. If $\{\Blq: 1\leq q\leq N(d,l)\}$ is an orthonormal basis of $\mathfrak{B}^l(\mbS^d)$, then it is known that $\{\Blq: l\in\mathbb{N}_0, 1\leq q\leq N(d,l)\}$ forms an orthonormal basis of $L^2((\mbS^d,\nu),\mathbb{C})$. Hereafter, $l$ that appears in $B^l_q$ or in other superscripts does not denote an exponent but denotes an index for notational simplicity. We note that the constant function $ B^0_1\equiv(\nu(\mbS^d))^{-1/2}$ is the orthonormal basis of $\mathfrak{B}^0(\mbS^d)$ since every spherical harmonic of order $0$ in $d+1$ variables is a constant function. Below, we summarize the examples of an orthonormal basis of $\mathfrak{B}^l(\mbS^d)$ for $l\in\mathbb{N}$.

\begin{example}\label{spherical harmonics example}\leavevmode
\begin{itemize}
\item[1.] ($d=1$) We note that each $x\in\mbS^1$ can be written as $x=(\cos\varphi_x,\sin\varphi_x)\tran$ for some $\varphi_x\in[0,2\pi)$. We define $\Blq:\mbS^1\rightarrow\mathbb{C}$ for $l\in\mathbb{N}$ by $ B^l_1(x)=\cos(l\varphi_x)/\sqrt{\pi}$ and $ B^l_2(x)=\sin(l\varphi_x)/\sqrt{\pi}$. Then, $\{\Blq: 1\leq q\leq2\}$ forms an orthonormal basis of $\mathfrak{B}^l(\mbS^1)$; see Chapter 2.2 in \cite{Atkinson and Han (2012)} for more details.
\item[2.] ($d=2$) We note that each $x\in\mbS^2$ can be written as
\begin{align*}
x=(\cos\varphi_x\sin\theta_x,\sin\varphi_x\sin\theta_x,\cos\theta_x)\tran
\end{align*}
for some $\varphi_x\in[0,2\pi)$ and $\theta_x\in[0,\pi)$. We define $\Blq:\mbS^2\rightarrow\mathbb{C}$ for $l\in\mathbb{N}$ by
\begin{align}\label{relation bar}
\Blq(x)=\sqrt{\frac{2l+1}{4\pi}}\cdot e^{\sqrt{-1}\cdot(q-l-1)\varphi_x}\cdot d^l_{q(l+1)}(\theta_x),
\end{align}
where $d^l_{qr}(\theta)\in\mathbb{R}$ for $1\leq q,r\leq2l+1$ and $\theta\in[0,\pi)$ is defined by
\begin{align}\label{small d matrix}
\begin{split}
    c^l_{qr}\cdot\sum_{k=\max\{0,r-q\}}^{\min\{2l+1-q,r-1\}}\frac{(-1)^{k+q-r}(\cos(\theta/2))^{2l-2k+r-q}(\sin(\theta/2))^{2k+q-r}}{(2l+1-q-k)!(r-1-k)!(k+q-r)!k!}
\end{split}
\end{align}
for $c^l_{qr}=((2l+1-q)!(q-1)!(2l+1-r)!(r-1)!)^{1/2}$. Then, $\{\Blq: 1\leq q\leq2l+1\}$ forms an orthonormal basis of $\mathfrak{B}^l(\mbS^2)$; see Theorem 2.1.1 in \cite{Terras (2013)}, Chapter 12.9 in \cite{Chirikjian (2012)} and Chapter 3.9 in \cite{Sakurai and Napolitano (2017)} for more details.
\item[3.] ($d\geq3$) An orthonormal basis of $\mathfrak{B}^l(\mbS^d)$ for $l\in\mathbb{N}$ and $d\geq3$ can be obtained recursively using an orthonormal basis of $\mathfrak{B}^{j}(\mbS^{d-1})$ for $0\leq j\leq l$. To describe this, we define the Legendre polynomial $P_{l,d+1}:[-1,1]\rightarrow\mathbb{R}$ of degree $l\in\mathbb{N}_0$ in $d+1$ variables by
    \[
    P_{l,d+1}(t)=l!\,\Gamma(d/2)\sum_{k=0}^{[l/2]}\frac{(-1)^k(1-t^2)^kt^{l-2k}}{4^kk!(l-2k)!\,\Gamma\left(k+\frac{d}{2}\right)}.
    \]
    We also define the normalized associated Legendre function $\tilde{P}_{l,d+1,j}:[-1,1]\rightarrow\mathbb{R}$ for $0\leq j\leq l$ by
    \[
    \tilde{P}_{l,d+1,j}(t)=\frac{((2l+d-1)(l+d+j-2)!)^{1/2}(1-t^2)^{j/2}}{2^{(d-1)/2+j}((l-j)!)^{1/2}\Gamma\left(j+\frac{d}{2}\right)}P_{l-j,d+1+2j}(t).
    \]
    We let $\{ B^{j}_{r}:1\leq r\leq N(j,d-1)\}$ be an orthonormal basis of $\mathfrak{B}^{j}(\mbS^{d-1})$ for $0\leq j\leq l$. We note that each $x\in\mbS^d$ can be written as
\begin{align}\label{spherical coordinate}
\begin{split}
x=\bigg(&\cos\varphi_x\prod_{k=1}^{d-1}\sin\theta_{kx},\sin\varphi_x\prod_{k=1}^{d-1}\sin\theta_{kx},\\
&\cos\theta_{1x}\prod_{k=2}^{d-1}\sin\theta_{kx},\cos\theta_{2x}\prod_{k=3}^{d-1}\sin\theta_{kx},\ldots,\cos\theta_{(d-1)x}\bigg)\tran
\end{split}
\end{align}
for some $\varphi_x\in[0,2\pi)$ and $\theta_{kx}\in[0,\pi)$ for $1\leq k\leq d-1$. We define $ B^l_{r,j}:\mbS^d\rightarrow\mathbb{C}$ for $l\in\mathbb{N}$ by
\begin{align*}
 B^l_{r,j}(x)=\tilde{P}_{l,d+1,j}(\cos\theta_{(d-1)x}) B^{j}_{r}(s(\varphi_x,\theta_{1x},\ldots,\theta_{(d-2)x})),
\end{align*}
where $s(\varphi_x,\theta_{1x},\ldots,\theta_{(d-2)x})$ is the point on $\mbS^{d-1}$ defined as the right hand side of (\ref{spherical coordinate}) with $d-1$ being replaced by $d-2$. Then, $\{ B^l_{r,j}:1\leq r\leq N(j,d-1), 0\leq j\leq l\}$ forms an orthonormal basis of $\mathfrak{B}^l(\mbS^d)$; see Chapter 2.11 in \cite{Atkinson and Han (2012)} for more details.
\end{itemize}
\end{example}

\subsection{Methodology}

In this section, we introduce our methodology. We let $X$ be a random vector taking values in $\mbS^d$ and $f_X$ be its density with respect to $\nu$. We assume that $f_X$ is square integrable. We let $SO(d+1)$ denote the space of all $(d+1)\times(d+1)$ real special orthogonal matrices. We recall that a $(d+1)\times(d+1)$ real matrix $A$ is called a special orthogonal matrix if $A\tran A=AA\tran=I_{d+1}$ and $\det(A)=1$, where $I_{d+1}$ is the $(d+1)\times(d+1)$ identity matrix. We suppose that we do not observe $X$ but we only observe $Z=UX$, where $U$ is an unobservable measurement error taking values in $SO(d+1)$ and $UX$ is the matrix multiplication between the matrix $U$ and vector $X$. We note that $Z\in\mbS^d$ since $\|UX\|^2=X\tran U\tran U X=X\tran X=1$. This measurement error is also natural since every matrix in $SO(d+1)$ rotates each point in $\mbS^d$ in a certain direction. For example, every matrix in $SO(2)$ can be written as
\begin{align}\label{so2 angle}
        \begin{psmallmatrix}
        \cos\varphi & -\sin\varphi\\
        \sin\varphi & \cos\varphi
        \end{psmallmatrix}
\end{align}
for some $\varphi\in[0,2\pi)$ and it rotates each point in $\mbS^1$ in the counter-clockwise direction by the angle $\varphi$. In addition, $SO(d+1)$ acts transitively on $\mbS^d$, i.e., for any $x_1, x_2\in\mbS^d$, there exists $u\in SO(d+1)$ such that $ux_1=x_2$. Our first aim is to estimate $f_X$ based on $n$ i.i.d. observations $\{Z_i:1\leq i\leq n\}$. Our second aim is to estimate the regression function $m:\mbS^d\rightarrow\mathbb{R}$ in the model
\begin{align}
Y=m(X)+\epsilon
\end{align}
based on $n$ i.i.d. observations $\{(Y_i,Z_i):1\leq i\leq n\}$, where $Y$ is a real-valued response and $\epsilon$ is an error term satisfying $\E(\epsilon|X)=0$. We note that this regression problem has not been covered for $d=2$ and $d\geq4$. Throughout this paper, we assume that $U$ is independent of $(X,\epsilon)$. This type of assumption is common in the literature of measurement error problems.

Below, we introduce a convolution property which is essential for our methodology. We define the convolution $g\ast f\in L^2((\mbS^d,\nu),\mathbb{C})$ of any two functions $g\in L^2((SO(d+1),\mu),\mathbb{C})$ and $f\in L^2((\mbS^d,\nu),\mathbb{C})$ by
\begin{align*}
(g\ast f)(x)=\int_{SO(d+1)}g(u)f(u^{-1}x)\,d\mu(u),
\end{align*}
where $u^{-1}$ is the inverse matrix of $u$, $u^{-1}x$ is the matrix multiplication between the matrix $u^{-1}$ and vector $x$, and $\mu$ is the normalized Haar measure on $SO(d+1)$. We recall that $\mu$ is the unique Borel probability measure on $SO(d+1)$ satisfying the left-translation-invariant property $\mu(A\mathcal{S})=\mu(\mathcal{S})$ for every $A\in SO(d+1)$ and Borel subset $\mathcal{S}\subset SO(d+1)$, where $A\mathcal{S}=\{AS:S\in\mathcal{S}\}$. We define the hyperspherical Fourier transform of $f$ at degree $l\in\mathbb{N}_0$ and order $1\leq q\leq N(d,l)$ by
\begin{align*}
\phi^l_q(f)=\int_{\mbS^d}f(x)\overline{\Blq(x)}d\nu(x).
\end{align*}
The hyperspherical Fourier transform is an analogue of the Euclidean Fourier transform with Euclidean domain, Lebesgue measure and Fourier basis function being replaced by $\mbS^d$, $\nu$ and $\Blq$, respectively. We also define a function $\mathscr{B}^l_q(\cdot;u):\mbS^d\rightarrow\mathbb{C}$ by $\mathscr{B}^l_q(x;u)=\Blq(ux)$ for $u\in SO(d+1)$. Since $\mathscr{B}^l_q(\cdot;u)$ belongs to $\mathfrak{B}^l(\mbS^d)$ (Proposition 4.7 in \cite{Efthimiou and Frye (2014)}) and $\{\Blq:1\leq q\leq N(d,l)\}$ is an orthonormal basis of $\mathfrak{B}^l(\mbS^d)$, it holds that
\begin{align}\label{relation sum}
\Blq(ux)=\sumr\lng\mathscr{B}^l_q(\cdot;u),\Blr\rng_2\Blr(x).
\end{align}
Finally, we define
\begin{align*}
\tilde{\phi}^l_{qr}(g)=\int_{SO(d+1)}g(u)D^l_{qr}(u)\,d\mu(u),
\end{align*}
where
\begin{align}\label{Dlqr}
D^l_{qr}(u)=\overline{\lng\mathscr{B}^l_q(\cdot;u),\Blr\rng_2}=\int_{\mbS^d}\overline{\Blq(ux)}\Blr(x)\,d\nu(x).
\end{align}
We call $\tilde{\phi}^l_{qr}(g)$ the $(q,r)$th element of the rotational Fourier transform of $g$ at degree $l$. Some practical examples of $\tilde{\phi}^l_{qr}(f_U)$ and $D^l_{qr}(u)$ are given in the Supplementary Material \ref{phiU implementation}. Then, the following convolution property holds.

\begin{proposition}\label{deconvolution general}
Let $f\in L^2((\mbS^d,\nu),\mathbb{C})$ and $g\in L^2((SO(d+1),\mu),\mathbb{C})$. Then, $\phi^l_q(g\ast f)=\sumr\tilde{\phi}^l_{qr}(g)\phi^l_r(f)$ for all $l\in\mathbb{N}_0$ and $1\leq q\leq N(d,l)$.
\end{proposition}

Proposition \ref{deconvolution general} is a generalization of Lemma 2.1 in \cite{Healy et al. (1998)} that considered the case where $d=2$. Now, we apply Proposition \ref{deconvolution general} to our setting. We let $f_U$ be the density of $U$ with respect to $\mu$. We assume that $f_U$ is square integrable. Then, one can show that the density $f_Z$ of $Z=UX\in\mbS^d$ with respect to $\nu$ exists and is given by $f_Z=f_U\ast f_X$. Hence, the following convolution property follows from Proposition \ref{deconvolution general}:
\begin{align}\label{deconvolution}
\phi^l_q(f_Z)=\sumr\tilde{\phi}^l_{qr}(f_U)\phi^l_r(f_X).
\end{align}
Defining $\phi^l(f)$ by the $N(d,l)$-vector whose $q$th element equals $\phi^l_q(f)$, and $\tilde{\phi}^l(g)$ by the $N(d,l)\times N(d,l)$ matrix whose $(q,r)$th element equals $\tilde{\phi}^l_{qr}(g)$, (\ref{deconvolution}) can be written as $\phi^l(f_Z)=\tilde{\phi}^l(f_U)\phi^l(f_X)$. Throughout this paper, we assume that the matrix $\tilde{\phi}^l(f_U)$ is invertible. Then, it can be rewritten as
\begin{align}\label{deconvolution matrix form}
\phi^l(f_X)=(\tilde{\phi}^l(f_U))^{-1}\phi^l(f_Z).
\end{align}
The invertibility of $\tilde{\phi}^l(f_U)$ is assumed in the literature of deconvolution density estimation on $\mbS^1$ or $\mbS^2$ (e.g. \cite{Efromovich (1997)}, \cite{Healy et al. (1998)}, \cite{Kim and Koo (2002)}, \cite{Kim et al. (2004)}). In fact, many popular distributions of $U$ such as the Laplace, Gaussian and von Mises-Fisher distributions on $SO(d+1)$ satisfy the invertibility. We give more concrete examples in the next section. In the literature of deconvolution density estimation on $\mbS^2$, it is always assumed that $f_U$ is known so that $\tilde{\phi}^l(f_U)$ is known. A Euclidean version of the latter assumption is also frequently assumed in the literature of Euclidean measurement error problems (e.g. \cite{Stefanski and Carroll (1990)}, \cite{Fan (1991a)}, \cite{Fan (1991b)}, \cite{Fan and Truong (1993)}, \cite{Delaigle et al. (2009)}, \cite{Belomestny and Goldenshluger (2021)}). Throughout this paper, we also focus on the case where $f_U$ is known, to build up a theoretical foundation in this new problem. This case is already challenging and the procedure starting from the case of known measurement error distribution has been adopted in the past new problems. In case $f_U$ is unknown, we may estimate it from additional data as in \cite{Efromovich (1997)}, \cite{Delaigle et al. (2008)}, \cite{Johannes (2009)}, \cite{Johannes and Schwarz (2013)} and \cite{Dattner et al. (2016)}, or by assuming a parametric distribution for $f_U$ and estimating its parameters without additional data as in \cite{Bertrand et al. (2019)}.

Now, we introduce our deconvolution estimator of $f_X$. Since $\{\Blq: l\in\mathbb{N}_0, 1\leq q\leq N(d,l)\}$ forms an orthonormal basis of $L^2((\mbS^d,\nu),\mathbb{C})$, it holds that
\begin{align}\label{Fourier series expansion original}
f_X=\sum_{l=0}^{\infty}\sumq\phi^l_q(f_X)\Blq
\end{align}
in the $L^2$ sense. The series at (\ref{Fourier series expansion original}) is called the Fourier-Laplace series of $f_X$. Under certain smoothness conditions on $f_X$, the series converges in the pointwise sense. We introduce such smoothness conditions in the next section. From (\ref{deconvolution matrix form}) and (\ref{Fourier series expansion original}), it holds that
\begin{align}\label{Fourier series expansion modified}
f_X=\sum_{l=0}^\infty\sumq\left(\sumr(\tilde{\phi}^l(f_U))_{qr}^{-1}\phi^l_r(f_Z)\right)\Blq,
\end{align}
where $(\tilde{\phi}^l(f_U))_{qr}^{-1}$ is the $(q,r)$th element of $(\tilde{\phi}^l(f_U))^{-1}$. Since $\phi^l_r(f_Z)=\E(\overline{\Blr(Z)})$ by definition, plugging the sample mean $n^{-1}\sum_{i=1}^n\overline{\Blr(Z_i)}$ in the place of $\phi^l_r(f_Z)$ at (\ref{Fourier series expansion modified}) gives an estimator of $f_X$. However, the estimator having the infinite sum $\sum_{l=0}^{\infty}$ is subject to a large variability since $(\tilde{\phi}^l(f_U))_{qr}^{-1}$ tends to infinity as $l\rightarrow\infty$. This tendency is analogous to the phenomenon in the Euclidean measurement error problems where the reciprocal of the Euclidean Fourier transform tends to infinity in the tails. To overcome this issue, we truncate the infinite sum $\sum_{l=0}^{\infty}$. Specifically, we let $0<T_n<\infty$ be a truncation level diverging to infinity as $n\rightarrow\infty$. Based on this truncation, we define
\begin{align}\label{deconvolution density estimator}
\hat{f}_X(x)=n^{-1}\sum_{i=1}^n\RE(K_{T_n}(x,Z_i)),
\end{align}
where
\begin{align*}
K_{T_n}(x,z)=\suml\sumq\left(\sumr(\tilde{\phi}^l(f_U))_{qr}^{-1}\overline{\Blr(z)}\right)\Blq(x)
\end{align*}
for $z\in\mbS^d$ and $\RE(K_{T_n}(x,Z_i))$ is the real part of $K_{T_n}(x,Z_i)$.
We note that $\hat{f}_X$ is the first density estimator in this general setting. \cite{Efromovich (1997)} and \cite{Healy et al. (1998)} introduced a similar density estimator defined by $n^{-1}\sum_{i=1}^nK_{T_n}(x,Z_i)$ for $d=1$ and $d=2$, respectively. However, $K_{T_n}(x,Z_i)$ is not necessarily real-valued, while $f_X$ is real-valued. Hence, it is natural to take $\RE(K_{T_n}(x,Z_i))$ instead of $K_{T_n}(x,Z_i)$ as in our density estimator $\hat{f}_X$. Taking the real part is also necessary to derive the asymptotic distribution of the estimator in Section \ref{rate and distribution} and the asymptotic confidence intervals for $f_X$ in Section \ref{confidence interval}. The below proposition tells that $\hat{f}_X$ is a reasonable density estimator in the sense that it integrates to one. This kind of property has not been noted for $d=2$ and $d\geq4$ in the literature.

\begin{proposition}\label{integral one}
$\int_{\mbS^d}K_{T_n}(x,z)d\nu(x)=1$ for all $z\in\mbS^d$, so that $\int_{\mbS^d}\RE(K_{T_n}(x,z))d\nu(x)=1$ for all $z\in\mbS^d$ and $\int_{\mbS^d}\hat{f}_X(x)d\nu(x)=1$.
\end{proposition}

We now introduce our deconvolution estimator of the regression function $m$. The proposed estimator of $m(x)=\E(Y|X=x)$ is given by
\begin{align}\label{deconvolution regression estimator}
\hat{m}(x)=\hat{f}_X(x)^{-1}n^{-1}\sum_{i=1}^n\RE(K_{T_n}(x,Z_i))Y_i.
\end{align}
We note that $\hat{m}$ is the first regression estimator for $d=2$ and $d\geq4$. Unlike the analysis of $\hat{f}_X(x)$, the analysis of $\hat{m}(x)$ requires an additional property on $\E(\RE(K_{T_n}(x,Z_i))|X_i)$ due to the additional term $Y_i$. To describe this, we define
\begin{align}\label{K star def}
K^*_{T_n}(x,x^*)=\suml\sumq\overline{\Blq(x^*)}\Blq(x)
\end{align}
for $x^*\in\mbS^d$. In view of (\ref{Fourier series expansion original}) and the fact $\E(\overline{\Blq(X)})=\phi^l_q(f_X)$, one may use $K^*_{T_n}$
to estimate $f_X$ in case the true values $\{X_i:1\leq i\leq n\}$ are observed. For instance, one may estimate $f_X(x)$ by
\begin{align*}
\hat{f}^*_X(x)=n^{-1}\sum_{i=1}^n\RE(K^*_{T_n}(x,X_i)),
\end{align*}
or simply by $n^{-1}\sum_{i=1}^nK^*_{T_n}(x,X_i)$, where the latter is the estimator studied by \cite{Hendriks (1990)} for the error-free case. One may also estimate $m(x)$ by
\begin{align}\label{naive regression estimator}
\hat{m}^*(x)=\hat{f}^*_X(x)^{-1}n^{-1}\sum_{i=1}^n\RE(K^*_{T_n}(x,X_i))Y_i.
\end{align}
The below proposition tells that $\E(K_{T_n}(x,Z)|X)$ equals $K^*_{T_n}(x,X)$. This kind of property has not been investigated for $d=2$ and $d\geq4$ in the literature.

\begin{proposition}\label{unbiased scoring}
$\E(K_{T_n}(x,Z)|X)=K^*_{T_n}(x,X)$ for all $x\in\mbS^d$, so that $\E(\RE(K_{T_n}(x,Z))|X)=\RE(K^*_{T_n}(x,X))$ for all $x\in\mbS^d$.
\end{proposition}

The above property that we term as `hyperspherical unbiased scoring' property gives that
\begin{align*}
\E\big(\hat m(x)\hat{f}_X(x)\big|X_1, \ldots, X_n\big)=\E\big(\hat m^*(x)\hat{f}^*_X(x)\big|X_1\ldots,X_n\big).
\end{align*}
The above identity tells that $K_{T_n}$ removes the effect of measurement errors in the bias of the nominator of $\hat{m}(x)$. We note that a Euclidean version of the above property was introduced in \cite{Stefanski and Carroll (1990)} and used in \cite{Fan and Truong (1993)} for the Euclidean measurement error problems. Such unbiased scoring properties are very important in regression analysis with measurement errors.

\section{Asymptotic properties}\label{rate and distribution}

\setcounter{equation}{0}
\setcounter{subsection}{0}

\subsection{Smoothness of measurement error distribution}\label{error distribution}

The asymptotic properties of our estimators depend on the smoothness of the measurement error density $f_U$. In the literature of Euclidean measurement error problems, two smoothness scenarios have been considered. They are ordinary-smooth and super-smooth scenarios; see \cite{Fan (1991a)}, for example. A typical example of ordinary-smooth distributions is the Laplace distribution, and a typical example of super-smooth distributions is the Gaussian distribution. In the literature of deconvolution density estimation on $\mbS^2$, three smoothness scenarios have been considered, namely ordinary-smooth, super-smooth and log-super-smooth scenarios. We extend the three scenarios to general $\mbS^d$. For this, we let $\|\cdot\|_\op$ denote the operator norm for complex matrices. For a $N(d,l)\times N(d,l)$ complex matrix $A$, it is defined by $\|A\|_{\op}=\sup\{\|Av\|_{\mathbb{C}^{N(d,l)}}:v\in\mathbb{C}^{N(d,l)}, \|v\|_{\mathbb{C}^{N(d,l)}}=1\}$, where $\|\cdot\|_{\mathbb{C}^{N(d,l)}}$ is the standard complex norm on $\mathbb{C}^{N(d,l)}$.

\begin{itemize}
\item[(S1)] (Ordinary-smooth scenario of order $\beta\geq0$) There exist constants $c_1, c_2>0$ such that, for all $l\in\mathbb{N}$, (i) $\|(\tilde{\phi}^l(f_U))^{-1}\|_\op\leq c_1\cdot l^\beta$ and (ii) $\|(\tilde{\phi}^l(f_U))^{-1}\|_\op\geq c_2\cdot l^{\beta}$.
\item[(S2)] (Super-smooth scenario of order $\beta>0$) There exist constants $c_1, c_2, \gamma>0$ and $\alpha\in\mathbb{R}$ such that, for all $l\in\mathbb{N}$, (i) $\|(\tilde{\phi}^l(f_U))^{-1}\|_\op\leq c_1\cdot l^\alpha\cdot\exp(\gamma\cdot l^\beta)$ and (ii) $\|(\tilde{\phi}^l(f_U))^{-1}\|_\op\geq c_2\cdot l^{\alpha}\cdot\exp(\gamma\cdot l^\beta)$.
\item[(S3)] (Log-super-smooth scenario of order $\beta>0$) There exist constants $c_1, c_2, \gamma>0$ and $\alpha, \xi_1, \xi_2\in\mathbb{R}$ such that, for all $l\in\mathbb{N}$, (i) $\|(\tilde{\phi}^l(f_U))^{-1}\|_\op\leq c_1\cdot l^\alpha\cdot\exp(\gamma\cdot l^\beta(\log{l}-\xi_1))$ and (ii) $\|(\tilde{\phi}^l(f_U))^{-1}\|_\op\geq c_2\cdot l^{\alpha}\cdot\exp(\gamma\cdot l^\beta(\log{l}-\xi_2))$.
\end{itemize}

We note that the conditions (i) in the scenarios (Sj) for $j\in\{1,2,3\}$ are for deriving the rates of convergence, while the conditions (ii) are used to verify a high-level condition for the asymptotic distributions. Distributions on $SO(d+1)$ are broadly studied in the literature (e.g. \cite{Leon et al. (2006)}, \cite{Sei et al. (2013)}, \cite{Qui et al. (2014)}, \cite{Nadarajah and Zhang (2017)}, \cite{Chakraborty and Vemuri (2019)}). We provide some examples of distributions satisfying the above scenarios. The ordinary-smooth scenario includes the case where there is no measurement error. In this case, $P(U=I_{N(d,l)})=1$, which gives $\tilde{\phi}^l(f_U)=I_{N(d,l)}$, where $I_{N(d,l)}$ is the $N(d,l)\times N(d,l)$ identity matrix. Hence, the case belongs to the ordinary-smooth scenario of order $\beta=0$. The Laplace distribution on $SO(d+1)$ with parameter $\lambda>0$ whose $\tilde{\phi}^l(f_U)$ is defined by $(1+\lambda^2\cdot l(l+d-1))^{-1}I_{N(d,l)}$ is another ordinary-smooth distribution. It satisfies (S1) with $\beta=2$. When $d=1$, its density is given by $f_U(u)=\pi(\exp(-\varphi_u/\lambda)/(1-\exp(-2\pi/\lambda))+\exp(\varphi_u/\lambda)/(\exp(2\pi/\lambda)-1))/\lambda$, where $\varphi_u\in[0,2\pi)$ is the angle corresponding to $u$ as given in (\ref{so2 angle}). When $d=2$, its density is given by $f_U(u)=\lambda^{-2}\pi\cos(a_\lambda(\pi-r_u))/(\cos(a_\lambda\pi)\sin(r_u/2))\cdot I(r_u>0)$, where $a_\lambda=\sqrt{1/4-\lambda^{-2}}\in\mathbb{C}$ and $r_u=\arccos((\Tr(u)-1)/2)\in[0,\pi]$ (Theorem 3.5 in \cite{Healy et al. (1998)}). Also, the Rosenthal distribution on $SO(3)$ with parameters $\theta\in(0,\pi]$ and $p>0$ whose density is given by $f_U(u)=\sum_{l=0}^\infty(2l+1)(\sin((2l+1)\theta/2)/((2l+1)\sin(\theta/2)))^p\sum_{q=-l}^lD^l_{qq}(u)$ has $\tilde{\phi}^l(f_U)=(\sin((2l+1)\theta/2)/((2l+1)\sin(\theta/2)))^pI_{2l+1}$ (\cite{Kim and Koo (2002)}), where $D^l_{qq}(u)$ is defined in (\ref{Dlqr}). Hence, it satisfies (S1) with $\beta=p$.

An example of super-smooth distributions is the Gaussian distribution on $SO(d+1)$ with parameter $\lambda>0$ whose $\tilde{\phi}^l(f_U)$ is defined by $\exp(-\lambda^2\cdot l(l+d-1)/2)I_{N(d,l)}$. It satisfies (S2) with $\beta=2, \alpha=0$ and $\gamma=\lambda^2/2$. When $d=1$, its density is given by $f_U(u)=\sqrt{2\pi}/\lambda\cdot\sum_{s\in\mathbb{Z}}\exp(-(\varphi_u+2\pi s)^2/(2\lambda^2))$. When $d=2$, its density is given by $f_U(u)=\sum_{l=0}^\infty(2l+1)\exp(-\lambda^2\cdot l(l+1)/2)\sum_{q=-l}^lD^l_{qq}(u)$.

Now, we consider the log-super-smooth scenario. Using Theorem 3 in \cite{Kim (2000)} or the result of Section 5.3 in \cite{Kim and Richards (2001)}, one may prove that the von Mises-Fisher distribution on $SO(d+1)$ with concentration parameter $\lambda>0$ and mean direction $A\in SO(d+1)$ is log-super-smooth. Its density is given by $f_U(u)=c(\lambda,A)^{-1}\exp(\lambda\cdot\Tr(A^{-1}u))$, where $c(\lambda,A)$ is the normalizing constant. When $d=1$, it satisfies (S3) with $\beta=1, \alpha=0, \gamma=1, \xi_1=1+\log\lambda$ and $\xi_2=1+\log(2\lambda)$. When $d=2$, it satisfies (S3) with $\beta=1, \alpha=4, \gamma=1, \xi_1=1+\log\lambda$ and $\xi_2=1+\log(3\lambda)$.

\subsection{Rates of convergence}\label{rates of convergence}

In this section, we discuss the uniform consistency and $L^2$ error rates of the density estimator $\hat{f}_X$ defined at (\ref{deconvolution density estimator}), and the $L^2$ error rates of the regression estimator $\hat{m}$ defined at (\ref{deconvolution regression estimator}). To state the required conditions, we denote the space of $s$-times continuously differentiable real-valued functions on $\mbS^d$ by $C^s(\mbS^d)$.

\begin{itemize}
\item[(A1)] For some $k\in\mathbb{N}$ with $k>d/4$, (i) $f_X\in C^{2k}(\mbS^d)$ and (ii) $m\in C^{2k}(\mbS^d)$.
\item[(A2)] (i) $f_X$ is bounded away from zero on $\mbS^d$ and (ii) $\E(Y^2|X=\cdot)$ is bounded on $\mbS^d$.
\end{itemize}

The condition (A1) is a smoothness condition on $f_X$ and $m$. Under (A1)-(i), the series at (\ref{Fourier series expansion original}) converges uniformly absolutely to $f_X$ by Theorem 2 in \cite{Kalf (1995)}. The uniform absolute convergence means that the absolute convergence holds uniformly. The condition (A2) is a standard regularity condition in nonparametric estimation. We also consider the following diverging speeds for the smoothing parameter $T_n$:

\begin{itemize}
\item[(T1)] (In the case of (S1)) $n^{-1/2}T_n^{\beta+d}=o(1)$.
\item[(T2)] (In the case of (S2)) $n^{-1/2}T_n^{\alpha+d}\exp(\gamma\cdot T_n^\beta)=o(1)$.
\item[(T3)] (In the case of (S3)) $n^{-1/2}T_n^{\alpha+d}\exp(\gamma\cdot T_n^\beta(\log T_n-\xi_1))=o(1)$.
\end{itemize}

Now, we are ready to state the asymptotic properties. We first introduce the uniform consistency of the density estimator $\hat{f}_X$. This is necessary to obtain the $L^2$ error rates of $\hat{m}$ and is also important in its own right.

\begin{proposition}\label{uniform consistency}
Assume that the Fourier-Laplace series at (\ref{Fourier series expansion original}) converges uniformly to $f_X$. Then, under either of the conditions (S1)-(i)+(T1), (S2)-(i)+(T2) and (S3)-(i)+(T3), it holds that
\begin{align*}
\sup_{x\in\mbS^d}|\hat{f}_X(x)-f_X(x)|=o_p(1).
\end{align*}
\end{proposition}

We note that the series at (\ref{Fourier series expansion original}) converges uniformly to $f_X$ under (A1)-(i) since the uniform absolute convergence implies the uniform convergence. Other weaker sufficient condition is that $f_X\in C^{s,\kappa}(\mbS^d)$ for some $s\geq0$ and $\kappa\in(0,1]$ with $s+\kappa>(d+1)/2-1$, where $C^{s,\kappa}(\mbS^d)$ is the space of real-valued functions whose $s$th order partial derivatives are H\"{o}lder continuous with exponent $\kappa$ (Theorem 2.36 in \cite{Atkinson and Han (2012)}). Now, we provide the $L^2$ rates of convergence for $\hat{f}_X$ and $\hat{m}$.

\begin{theorem}\label{L2 rates}
Assume that the condition (A1)-(i) holds. Then,
\begin{itemize}
\item[(a)] Under (S1)-(i) and (T1), it holds that
\begin{align*}
\int_{\mbS^d}|\hat{f}_X(x)-f_X(x)|^2d\nu(x)=O_p(T_n^{-4k}+n^{-1}T_n^{2\beta+d}).
\end{align*}
The same rate holds for $\int_{\mbS^d}|\hat{m}(x)-m(x)|^2d\nu(x)$ under the additional conditions (A1)-(ii) and (A2).
\item[(b)] Under (S2)-(i) and (T2), it holds that
\begin{align*}
\int_{\mbS^d}|\hat{f}_X(x)-f_X(x)|^2d\nu(x)=O_p(T_n^{-4k}+n^{-1}T_n^{2\alpha+d}\exp(2\gamma\cdot T_n^\beta)).
\end{align*}
The same rate holds for $\int_{\mbS^d}|\hat{m}(x)-m(x)|^2d\nu(x)$ under the additional conditions (A1)-(ii) and (A2).
\item[(c)] Under (S3)-(i) and (T3), it holds that
\begin{align*}
\int_{\mbS^d}|\hat{f}_X(x)-f_X(x)|^2d\nu(x)=O_p(T_n^{-4k}+n^{-1}T_n^{2\alpha+d}\exp(2\gamma\cdot T_n^\beta(\log{T_n}-\xi_1))).
\end{align*}
The same rate holds for $\int_{\mbS^d}|\hat{m}(x)-m(x)|^2d\nu(x)$ under the additional conditions (A1)-(ii) and (A2).
\end{itemize}
\end{theorem}

We note that the above $L^2$ error rates converge to zero as $n\rightarrow\infty$. The term $T_n^{-4k}$ in the rates comes from the bias parts of the estimators, and the remaining term in each rate is originated from a stochastic part contributing to the variance. We may optimize each $L^2$ error rate by taking a suitable speed of $T_n\rightarrow\infty$. Specifically, we consider the following speeds:

\begin{itemize}
\item[(T1$'$)] (In the case of (S1)) $T_n\asymp n^{1/(4k+2\beta+d)}$.
\item[(T2$'$)] (In the case of (S2)) $T_n=K\cdot(\log{n})^{1/\beta}$ for $0<K<(2\gamma)^{-1/\beta}$.
\item[(T3$'$)] (In the case of (S3)) $T_n=K\cdot(\log{n}/\log{\log{n}})^{1/\beta}$ for $0<K<(2\gamma/\beta)^{-1/\beta}$.
\end{itemize}

The speed (T1$'$) is optimal in the sense that it balances the asymptotic bias $T_n^{-4k}$ and asymptotic variance $n^{-1}T_n^{2\beta+d}$. In the cases of super-smoothness and log-super-smoothness, however, there exists no such speed that makes the corresponding asymptotic bias and variance be of the same magnitude. This is because $T_n$ also appears in the exponents $\exp(2\gamma\cdot T_n^\beta)$ and $\exp(2\gamma\cdot T_n^\beta(\log{T_n}-\xi_1))$, respectively. The choices of $T_n$ given in (T2$'$) and (T3$'$) have specific constant factors $K$ with constraints. The upper bounds of $K$ are actually the thresholds, beyond which $n^{-1}T_n^{2\alpha+d}\exp(2\gamma\cdot T_n^\beta)$
and $n^{-1}T_n^{2\alpha+d}\exp(2\gamma\cdot T_n^\beta(\log{T_n}-\xi_1))$, respectively, diverge to infinity, while they are dominated
by $T_n^{-4k}$ for $K$ smaller than the thresholds. We note that similar constraints have been put on bandwidths in the Euclidean super-smooth scenario; see Theorem 1 in \cite{Fan (1991a)}, and Theorem 1 and Remark 1 in \cite{Fan and Truong (1993)}, for example.

\begin{corollary}\label{L2 rates optimal}
Assume that the condition (A1)-(i) holds. Then,
\begin{itemize}
\item[(a)] Under (S1)-(i) and (T1$\,'$), it holds that
\begin{align*}
\int_{\mbS^d}|\hat{f}_X(x)-f_X(x)|^2d\nu(x)=O_p(n^{-4k/(4k+2\beta+d)}).
\end{align*}
The same rate holds for $\int_{\mbS^d}|\hat{m}(x)-m(x)|^2d\nu(x)$ under the additional conditions (A1)-(ii) and (A2).
\item[(b)] Under (S2)-(i) and (T2$\,'$), it holds that
\begin{align*}
\int_{\mbS^d}|\hat{f}_X(x)-f_X(x)|^2d\nu(x)=O_p((\log{n})^{-4k/\beta}).
\end{align*}
The same rate holds for $\int_{\mbS^d}|\hat{m}(x)-m(x)|^2d\nu(x)$ under the additional conditions (A1)-(ii) and (A2).
\item[(c)] Under (S3)-(i) and (T3$\,'$), it holds that
\begin{align*}
\int_{\mbS^d}|\hat{f}_X(x)-f_X(x)|^2d\nu(x)=O_p((\log{n}/\log{\log{n}})^{-4k/\beta}).
\end{align*}
The same rate holds for $\int_{\mbS^d}|\hat{m}(x)-m(x)|^2d\nu(x)$ under the additional conditions (A1)-(ii) and (A2).
\end{itemize}
\end{corollary}

It is natural that the rate in (a) of Corollary \ref{L2 rates optimal} gets slower as the dimension $d$ increases, due to the well known phenomena called the curse of dimensionality. However, the rates in (b) and (c) of Corollary \ref{L2 rates optimal} are independent of $d$ since the rates are dominated by the asymptotic bias $T_n^{-4k}$ which is independent of $d$.
We note that similar log-type error rates were obtained by \cite{Fan (1991a)} and \cite{Fan and Truong (1993)} for the Euclidean measurement error problems under the Euclidean super-smooth scenario.

\subsection{Asymptotic distributions}\label{asymptotic distributions}

In this section, we discuss the asymptotic distributions of our density and regression estimators. Recently, \cite{Jeon et al. (2021)} derived some asymptotic distributions for their deconvolution estimators on compact and connected Lie groups. We note that $\mbS^1$ and $\mbS^3$ are such Lie groups. To the best of our knowledge, however, no asymptotic distribution has been derived for $\mbS^d$ with $d=2$ and $d\geq4$ in the literature of measurement error problems.

We first derive the asymptotic distribution of $\hat{f}_X$. Before we state the result, we introduce a high-level condition. In the following high level condition, $a_n\gtrsim b_n$ for two positive sequences $a_n$ and $b_n$ means that there exists a constant $c>0$ such that $a_n\geq c\cdot b_n$ for all $n$. We also define $a_n\lesssim b_n$ in the obvious way.

\begin{itemize}
\item[(B1)] (In the case of (S1)) There exists a constant $0\leq q\leq d$ such that, for each $x\in\mbS^d$, $\E((\RE(K_{T_n}(x,Z)))^2)\gtrsim T_n^{2\beta+q}$.
\item[(B2)] (In the case of (S2)) There exists a constant $0\leq q\leq d$ such that, for each $x\in\mbS^d$ and $0<\eta<1$, $\E((\RE(K_{T_n}(x,Z)))^2)\gtrsim T_n^{2\alpha+q}\exp(2\gamma(\eta\cdot T_n)^\beta)$.
\item[(B3)] (In the case of (S3)) There exist constants $0\leq q\leq d$ and $\zeta\in\mathbb{R}$ such that, for each $x\in\mbS^d$ and $0<\eta<1$, $\E((\RE(K_{T_n}(x,Z)))^2)\gtrsim T_n^{2\alpha+q}\exp(2\gamma(\eta\cdot T_n)^\beta(\log T_n-\zeta))$.
\end{itemize}

The lower bounds to $\E\big(\RE(K_{T_n}(x,Z))^2\big)$ in the conditions (B1)-(B3) are motivated by the upper bounds to $\E((\RE(K_{T_n}(x,Z)))^2)$, which are of the magnitude
\[
T_n^{2\beta+d}, \quad T_n^{2\alpha+d}\exp(2\gamma\cdot T_n^\beta) \quad \mbox{or} \quad T_n^{2\alpha+d}\exp(2\gamma\cdot T_n^\beta(\log T_n-\zeta))
\]
depending on the smoothness scenarios. The verification of (B1)-(B3)
with $q=d$ is particularly important in constructing asymptotic confidence intervals for $f_X$ and $m$. We verify them for $d\in\{1,2\}$ and various $f_U$ in the next section. We put the range $0\leq q\leq d$ in (B1)-(B3) to give flexibility. We also give more flexibility in choosing $T_n$ instead of choosing the specific ones in (T1$'$)-(T3$'$). The following flexible speeds cover the speeds in (T1$'$)-(T3$'$).

\begin{itemize}
\item[(T1$''$)] $T_n\asymp n^p$ for some $0<p<1/(2d-q)$, where $q$ is the constant in (B1).
\item[(T2$''$)] $T_n\lesssim (\log n)^{1/\beta}$ for $\beta$ in (S2).
\item[(T3$''$)] $T_n\lesssim (\log n/\log \log n)^{1/\beta}$ for $\beta$ in (S3).
\end{itemize}

\begin{theorem}\label{asymptotic distribution density}
Assume that the Fourier-Laplace series at (\ref{Fourier series expansion original}) converges pointwise to $f_X$. Then, under either of the conditions (S1)-(i)+(T1$\,''$)+(B1), (S2)-(i)+(T2$\,''$)+(B2) and (S3)-(i)+(T3$\,''$)+(B3), it holds that, for all $x\in \mbS^d$,
\begin{align*}
\sqrt{n}\cdot\frac{\hat{f}_X(x)-f_X(x)-\big(\E(\RE(K_{T_n}(x,Z)))-f_X(x)\big)}{\sqrt{\Var(\RE(K_{T_n}(x,Z)))}}\overset{d}{\longrightarrow}N(0,1).
\end{align*}
\end{theorem}

We note that the condition on the Fourier-Laplace series in Theorem \ref{asymptotic distribution density} is weaker than the corresponding condition in Proposition \ref{uniform consistency}. Now, we investigate the asymptotic distribution of $\hat{m}$ in the ordinary-smooth scenario. Deriving it for the super-smooth and log-super-smooth scenarios has a technical issue. The issue is that
\begin{align}\label{normality difficulty}
\sqrt{n}\cdot\frac{\E(\RE(K_{T_n}(x,Z))(Y-m(x)))\cdot(\hat{f}_X(x)-f_X(x))}{\sqrt{\E((\RE(K_{T_n}(x,Z)))^2)}}=o_p(1)
\end{align}
does not hold in the super-smooth and log-super-smooth scenarios. (\ref{normality difficulty}) is an important part in the proof; see Remark \ref{only ordinary} immediately after the proof of Theorem \ref{asymptotic distribution regression} for more details. For the asymptotic distribution of $\hat{m}$ in the ordinary-smooth scenario, we make an additional condition.

\begin{itemize}
\item[(B4)] $\E(|Y|^{2+\delta}|X=\cdot)$ is bounded on $\mbS^d$ for some $\delta>0$ and $\E(\epsilon^2|X=\cdot)$ is bounded away from zero on $\mbS^d$.
\end{itemize}

The condition (B4) is a standard regularity condition in nonparametric regression. We also consider a new flexible range on $T_n$ for the ordinary smooth scenario. The following range based on the constants $\beta$ in (S1) and $k$ in (A1) covers the speed in (T1$'$).

\begin{itemize}
\item[(T1$'''$)] $T_n\asymp n^p$ for some $1/(2\beta+d+8k)\leq p<1/(2\beta+2d)$.
\end{itemize}

We note that the range of $p$ in (T1$'''$) is valid since $k$ in (A1) satisfies $k>d/4$. The upper bound $1/(2\beta+2d)$ in the range is required to make $T_n$ satisfy (T1). To state the next theorem, we denote by $\Delta_{\mbS^d}$ the Laplace-Beltrami operator associated with $\mbS^d$ (twice differential operator acting on twice continuously differentiable functions on $\mbS^d$), and by $\Delta_{\mbS^d}^s$ the compositions of $\Delta_{\mbS^d}$ for $s$ times. We note that, if a function $f:\mbS^d\rightarrow\mathbb{R}$ is $2s$-times continuously differentiable on $\mbS^d$, then the function $\Delta_{\mbS^d}^s(f):\mbS^d\rightarrow\mathbb{R}$ is well defined and is continuous on $\mbS^d$.

\begin{theorem}\label{asymptotic distribution regression}
Assume that the conditions (S1)-(i), (A1), (A2)-(i), (B1) with $q=d$, (B4) and (T1$\,'''$) hold, and that the Fourier-Laplace series of $\Delta_{\mbS^d}^k(f_X)$ and of $\Delta_{\mbS^d}^k(m\cdot f_X)$ converge absolutely on $\mbS^d$. Then, it holds that, for all $x\in \mbS^d$,
\begin{align*}
\sqrt{n}\cdot\frac{\hat{m}(x)-m(x)-\E(\RE(K_{T_n}(x,Z))(Y-m(x)))/f_X(x)}{\sqrt{\Var(\RE(K_{T_n}(x,Z))(Y-m(x)))}/f_X(x)}\overset{d}{\longrightarrow}N(0,1).
\end{align*}
\end{theorem}

Theorem \ref{asymptotic distribution regression} also holds for general $q$ in (B1) with more complex versions of (A1) and (T1$'''$) although we only state it with $q=d$ for simplicity. Regarding the condition on the absolute convergence in Theorem \ref{asymptotic distribution regression}, we recall that, if a function $f:\mbS^d\rightarrow\mathbb{R}$ is $2s$-times continuously differentiable for $s>d/4$, then its Fourier-Laplace series is absolutely convergent. However, for certain $d$, much weaker sufficient conditions exist. For example, the Fourier-Laplace series of a function on $\mbS^1$ is absolutely convergent if the function is H\"{o}lder continuous with exponent greater than 1/2 or if the function is of bounded variation and H\"{o}lder continuous with positive exponent; see \cite{Katznelson (2004)}.

\section{Asymptotic confidence intervals}\label{confidence interval}

\setcounter{equation}{0}
\setcounter{subsection}{0}

In this section, we verify the high-level conditions (B1)-(B3) with $q=d$ for certain cases and provide two types of asymptotic confidence intervals for both $f_X$ and $m$. One type is based on the asymptotic normality given in Theorems~\ref{asymptotic distribution density} and \ref{asymptotic distribution regression},
and the other is based on empirical likelihoods. For the first type, we estimate the biases $\E\big(\RE(K_{T_n}(x,Z))\big)-f_X(x)$ in the nominator in Theorem \ref{asymptotic distribution density} and $\E\big(\RE\big(K_{T_n}(x,Z)\big)\big(Y-m(x)\big)\big)/f_X(x)$ in the nominator in Theorem \ref{asymptotic distribution regression}.
We estimate them simply by zero. These are natural choices since plugging $\hat f_X(x)=n^{-1}\sum_{i=1}^n \RE(K_{T_n}(x,Z_i))$
into both $\E\big(\RE(K_{T_n}(x,Z))\big)$ and $f_X(x)$ gives zero estimates for $\E\big(\RE(K_{T_n}(x,Z))\big)-f_X(x)$ and plugging $0=n^{-1}\sum_{i=1}^n \RE\big(K_{T_n}(x,Z_i)\big)\big(Y_i-\hat m(x)\big)$ into
$\E\big(\RE\big(K_{T_n}(x,Z)\big)\big(Y-m(x)\big)\big)$ gives zero estimates for $\E\big(\RE\big(K_{T_n}(x,Z)\big)\big(Y-m(x)\big)\big)/f_X(x)$.

To justify the zero estimates of the biases in the construction of the first type asymptotic confidence intervals, it is essential to verify that the variances dominate the squared biases.
In the case of $\hat f_X$, this amounts to showing
\begin{align}\label{negligible2}
\sqrt{n}\cdot\frac{\E\big(\RE(K_{T_n}(x,Z))\big)-f_X(x)}{\sqrt{\E\big(\big(\RE(K_{T_n}(x,Z))\big)^2\big)}}=o(1)
\end{align}
since $\E\big(\big(\RE(K_{T_n}(x,Z))\big)^2\big)$ determines the magnitude of $\Var(\RE(K_{T_n}(x,Z)))$; see the proof of Theorem \ref{asymptotic distribution density toruses}. For the verification of (\ref{negligible2}), we need sharp lower bounds to the denominator. It can be accomplished by verifying (B1)-(B3) with $q=d$.

\subsection{Verification of (B1)-(B3)}

In this section, we verify that $\E((\RE(K_{T_n}(x,Z)))^2)$ achieves the lower bounds given in (B1)-(B3) with $q=d$. In the Euclidean nonparametric statistics, a common approach to get a lower bound of such quantity is to find an asymptotic leading term by applying Taylor expansion. However, since there exists no suitable Taylor expansion for our problem, it is not trivial to verify (B1)-(B3) not only with the maximal $q=d$ but also with the minimal $q=0$.

We first show that $\E(|K_{T_n}(x,Z)|^2)$ achieves the lower bounds in (B1)-(B3) with $q=d$, and then consider the lower bounds to $\E((\RE(K_{T_n}(x,Z))^2))$. For this, we need an additional condition. We let $\sigma_{\rm min}((\tilde{\phi}^l(f_U))^{-1})$ denote the minimum singular value of $(\tilde{\phi}^l(f_U))^{-1}$.

\begin{itemize}
\item[(C)] There exists a positive constant $c$ such that, for all $l\in\mathbb{N}_0$, $\|(\tilde{\phi}^l(f_U))^{-1}\|_{\rm op}\leq c\cdot\sigma_{\rm min}((\tilde{\phi}^l(f_U))^{-1})$.
\end{itemize}

One can easily check that the condition (C) holds with $c=1$ for any $f_U$ on $\mbS^1$. It also holds with $c=1$ for any $f_U$ satisfying $\tilde{\phi}^l(f_U)=s_l\cdot I_{N(d,l)}$ for some $0\neq s_l\in\mathbb{R}$. Examples of such distributions for $d\geq2$ include the Laplace, Rosenthal, Gaussian and error-free distributions that we introduced immediately after (S1)-(S3).

\begin{lemma}\label{rate of denominator complex}
Assume that the conditions (A2)-(i) and (C) hold. Then, for each $j\in\{1,2,3\}$, under (Sj)-(ii),
$\E(|K_{T_n}(x,Z)|^2)$ attains the lower bound given in (Bj) with $q=d$.
\end{lemma}

Lemma \ref{rate of denominator complex} tells about the lower bounds to $\E(|K_{T_n}(x,Z)|^2)$. However, since
\begin{align*}
\E(|K_{T_n}(x,Z)|^2)=\E((\RE(K_{T_n}(x,Z)))^2)+\E((\IM(K_{T_n}(x,Z)))^2)\geq\E((\RE(K_{T_n}(x,Z))^2)),
\end{align*}
where $\IM(K_{T_n}(x,Z))$ is the imaginary part of $K_{T_n}(x,Z)$, Lemma \ref{rate of denominator complex} does not provide the lower bounds to $\E((\RE(K_{T_n}(x,Z))^2))$. This causes another difficulty. Our original attempt for this issue was to show that $\E((\RE(K_{T_n}(x,Z)))^2)\geq\E((\IM(K_{T_n}(x,Z)))^2)$ always holds, but it was not successful due to computational difficulties. Instead, we managed to prove that
\begin{align}\label{integral comparison}
\int_{\mbS^d}\RE(K_{T_n}(x,z))^2d\nu(z)\geq\int_{\mbS^d}\IM(K_{T_n}(x,z))^2d\nu(z)
\end{align}
for certain cases. Since $\int_{\mbS^d}|K_{T_n}(x,z)|^2d\nu(z)$ also achieves the same lower bounds as those to $\E(|K_{T_n}(x,Z)|^2)$ as demonstrated in the proof of Lemma \ref{rate of denominator complex}, proving (\ref{integral comparison}) provides that $\int_{\mbS^d}\RE(K_{T_n}(x,z))^2d\nu(z)$ also achieves the same lower bounds. This with the assumption $\inf_{z\in\mbS^d}f_Z(z)>0$ gives the desired lower bounds to $\E((\RE(K_{T_n}(x,Z)))^2)$. We note that the latter assumption on the density $f_Z$ of $Z$ is implied by (A2)-(i). The cases for which (\ref{integral comparison}) hold are the followings:

\begin{itemize}
\item[(G1)] $d=1$.
\item[(G2)] $d=2$ and $\tilde{\phi}^l(f_U)=s_l\cdot I_{N(d,l)}$ for some $0\neq s_l\in\mathbb{R}$ and for all $l\in\mathbb{N}_0$.
\end{itemize}

Verification of (\ref{integral comparison}) for the cases (G1)-(G2) requires complex computation based on the theory of spherical harmonics. We note that (G1)-(G2) are practically the most important cases. They cover circular data and spherical data. The distributions of $U$ satisfying (G2) include, but are not limited to, the Laplace, Rosenthal, Gaussian and error-free distributions on $SO(3)$. We also note that (G1)-(G2) satisfy the condition (C). Hence, (B1)-(B3) with $q=d$ follow for the cases (G1)-(G2) under the condition (A2)-(i).

\begin{lemma}\label{rate of denominator toruses}
Assume that the conditions (A2)-(i) and either (G1) or (G2) hold. Then, for each $j\in\{1,2,3\}$, under (Sj)-(ii), (Bj) holds with $q=d$.
\end{lemma}

Although we verify (B1)-(B3) with $q=d$ only for the above cases due to computational difficulties, we strongly believe that they hold for general $\mbS^d$ and $f_U$. For the verification, one could apply a special computation technique or a different way without direct computation that we are currently not aware. We leave it as an open problem.

\subsection{Confidence intervals based on asymptotic normality}

In this section, we construct asymptotic confidence intervals based on Theorems \ref{asymptotic distribution density} and \ref{asymptotic distribution regression} for the ordinary-smooth scenario under the condition (B1) with $q=d$. We only treat the ordinary-smooth scenario since (\ref{negligible2}) does not hold with the speeds of $T_n$ in (T2$''$) and (T3$''$) in the super-smooth and log-super-smooth scenarios; see Remark \ref{only ordinary} in the Supplementary Material for a related discussion. Even in the Euclidean measurement error problems, studying the Euclidean ordinary-smooth scenario is more common than studying the Euclidean super-smooth scenario due to many technical difficulties in the Euclidean super-smooth scenario.

For the construction of the asymptotic confidence intervals, we need to estimate the unknown biases and variances in Theorems \ref{asymptotic distribution density}-\ref{asymptotic distribution regression}. The biases are $\E(\RE(K_{T_n}(x,Z)))-f_X(x)$ and $\E(\RE(K_{T_n}(x,Z))(Y-m(x)))/f_X(x)$ and variances are $s^2_1(x)$ and $s^2_2(x)$, where
\begin{align*}
s_1(x)&=\big(\Var(\RE(K_{T_n}(x,Z)))\big)^{1/2},\\
s_2(x)&=\big(\Var(\RE(K_{T_n}(x,Z))(Y-m(x)))\big)^{1/2}/f_X(x).
\end{align*}
We estimate the biases by zero as we demonstrated in the beginning of Section \ref{confidence interval}. For the variances, we use the following natural estimators:
\begin{align*}
\hat{s}^2_1(x)=&n^{-1}\sum_{i=1}^n(\RE(K_{T_n}(x,Z_i)))^2-(\hat{f}_X(x))^2,\\
\hat{s}^2_2(x)=&\hat{f}^{-2}_X(x)\cdot\bigg(n^{-1}\sum_{i=1}^n(\RE(K_{T_n}(x,Z_i))(Y_i-\hat{m}(x)))^2\\
&\qquad\qquad\qquad-\bigg(n^{-1}\sum_{i=1}^n\RE(K_{T_n}(x,Z_i))(Y_i-\hat{m}(x))\bigg)^2\bigg)\\
=&\hat{f}^{-2}_X(x)\cdot n^{-1}\sum_{i=1}^n(\RE(K_{T_n}(x,Z_i))(Y_i-\hat{m}(x)))^2.
\end{align*}

\begin{theorem}\label{asymptotic distribution density toruses}
Assume that the conditions (S1)-(i), (A1)-(i), (T1$\,'$) and (B1) with $q=d$ hold, and that the Fourier-Laplace series of $\Delta_{\mbS^d}^k(f_X)$ converges absolutely on $\mbS^d$. Then, it holds that, for all $x\in\mbS^d$,
\begin{align*}
\sqrt{n}\cdot\frac{\hat{f}_X(x)-f_X(x)}{\hat{s}_1(x)}\overset{d}{\longrightarrow}N(0,1).
\end{align*}
Hence, a $(1-\alpha)\times100\%$ asymptotic confidence interval for $f_X(x)$ is given by
\begin{align*}
\left(\hat{f}_X(x)-z_{\alpha/2}\frac{\hat{s}_1(x)}{\sqrt{n}},\hat{f}_X(x)+z_{\alpha/2}\frac{\hat{s}_1(x)}{\sqrt{n}}\right).
\end{align*}
\end{theorem}

\begin{theorem}\label{asymptotic distribution regression toruses}
Assume that the conditions (S1)-(i), (A1), (A2)-(i), (T1$\,'$), (B4) and (B1) with $q=d$ hold, and that the Fourier-Laplace series of $\Delta_{\mbS^d}^k(f_X)$ and of $\Delta_{\mbS^d}^k(m\cdot f_X)$ converge absolutely on $\mbS^d$. Then, it holds that, for all $x\in\mbS^d$,
\begin{align*}
\sqrt{n}\cdot\frac{\hat{m}(x)-m(x)}{\hat{s}_2(x)}\overset{d}{\longrightarrow}N(0,1).
\end{align*}
Hence, a $(1-\alpha)\times100\%$ asymptotic confidence interval for $m(x)$ is given by
\begin{align*}
\left(\hat{m}(x)-z_{\alpha/2}\frac{\hat{s}_2(x)}{\sqrt{n}},\hat{m}(x)+z_{\alpha/2}\frac{\hat{s}_2(x)}{\sqrt{n}}\right).
\end{align*}
\end{theorem}

\subsection{Confidence intervals based on empirical likelihoods}

Asymptotic confidence regions based on empirical likelihoods, called empirical likelihood confidence regions, are
useful alternatives to those based on asymptotic normality. Empirical likelihood confidence regions have the advantage
that their shape is determined by the data, they are invariant under transformations,
and they often do not require the estimation of the variance. For an introduction to empirical likelihood methods, we refer to \cite{Owen (2001)}. A broad review and a general theory for empirical likelihood methods can be found in \cite{Chen and Van Keilegom (2009)} and \cite{Hjort et al. (2009)}, respectively. To construct the empirical likelihood confidence regions for $f_X$ and $m$, we define $F_{f_X}(Z_i,\theta;x)=\RE(K_{T_n}(x,Z_i))-\theta$ and $F_m(Z_i,Y_i,\theta;x)=\RE(K_{T_n}(x,Z_i))(Y_i-\theta)$ for $\theta\in\mathbb{R}$ and $x\in\mbS^d$. We also define the corresponding empirical likelihood ratio functions at $x$ by
\begin{align*}
\EL_{f_X}(\theta;x)&=\max\left\{\prod_{i=1}^n(nw_i):w_i>0, \sum_{i=1}^nw_i=1, \sum_{i=1}^nw_iF_{f_X}(Z_i,\theta;x)=0\right\},\\
\EL_m(\theta;x)&=\max\left\{\prod_{i=1}^n(nw_i):w_i>0, \sum_{i=1}^nw_i=1, \sum_{i=1}^nw_iF_m(Z_i,Y_i,\theta;x)=0\right\}.
\end{align*}
Here, we define the maximum of the empty set to be zero. Then, we define the respective empirical likelihood confidence regions by $\{\theta\in\mathbb{R}:\EL_{f_X}(\theta;x)\geq c_{f_X}\}$ and $\{\theta\in\mathbb{R}:\EL_m(\theta;x)\geq c_m\}$ for some positive constants $c_{f_X}$ and $c_m$. To determine the constants, we provide the asymptotic distributions of the empirical likelihood ratio functions. For this, we take the following conditions.

\begin{itemize}
\item[(E1)] $P(\EL_{f_X}(f_X(x);x)>0)\rightarrow1$ for each $x\in\mbS^d$.
\item[(E2)] $P(\EL_m(m(x);x)>0)\rightarrow1$ for each $x\in\mbS^d$.
\end{itemize}

The conditions (E1)-(E2) are basic in the empirical likelihood technique. We note that $\EL_{f_X}(f_X(x);x)>0$ is satisfied as long as there are at least two data points $Z_i$ and $Z_j$ such that $F_{f_X}(Z_i,f_X(x);x)>0$ and $F_{f_X}(Z_j,f_X(x);x)<0$, and $\EL_m(m(x);x)>0$ is satisfied as long as there are at least two data points $(Z_i,Y_i)$ and $(Z_j,Y_j)$ such that $F_m(Z_i,Y_i,m(x);x)>0$ and $F_m(Z_j,Y_j,m(x);x)<0$. We also treat the ordinary-smooth scenario only since similar technical difficulties exist. Below, $\chi^2_\alpha(1)$ denotes the $(1-\alpha)$ quantile of the chi-square distribution $\chi^2(1)$ of degree 1.

\begin{theorem}\label{empirical likelihood density}
Assume that the conditions (S1)-(i), (A1)-(i), (T1$\,'$), (E1) and (B1) with $q=d$ hold, and that the Fourier-Laplace series of $\Delta_{\mbS^d}^k(f_X)$ converges absolutely on $\mbS^d$. Then, it holds that, for all $x\in\mbS^d$,
\begin{align*}
-2\log\EL_{f_X}(f_X(x);x)\overset{d}{\longrightarrow}\chi^2(1).
\end{align*}
Hence, a $(1-\alpha)\times100\%$ asymptotic confidence region for $f_X(x)$ is given by $\{\theta\in\mathbb{R}:-2\log\EL_{f_X}(\theta;x)\leq\chi^2_\alpha(1)\}=\{\theta\in\mathbb{R}:\EL_{f_X}(\theta;x)\geq\exp(-\chi^2_\alpha(1)/2)\}$.
\end{theorem}

\begin{theorem}\label{empirical likelihood regression}
Assume that the conditions (S1)-(i), (A1), (A2)-(i), (T1$\,'$), (B4), (E2) and (B1) with $q=d$ hold, and that the Fourier-Laplace series of $\Delta_{\mbS^d}^k(f_X)$ and of $\Delta_{\mbS^d}^k(m\cdot f_X)$ converge absolutely on $\mbS^d$. Then, it holds that, for all $x\in\mbS^d$,
\begin{align*}
-2\log\EL_m(m(x);x)\overset{d}{\longrightarrow}\chi^2(1).
\end{align*}
Hence, a $(1-\alpha)\times100\%$ asymptotic confidence region for $m(x)$ is given by $\{\theta\in\mathbb{R}:-2\log\EL_m(\theta;x)\leq\chi^2_\alpha(1)\}=\{\theta\in\mathbb{R}:\EL_m(\theta;x)\geq\exp(-\chi^2_\alpha(1)/2)\}$.
\end{theorem}

We note that the asymptotic confidence regions in Theorems \ref{empirical likelihood density} and \ref{empirical likelihood regression} are in fact intervals. This is because $t\theta_1+(1-t)\theta_2$ for $0<t<1$ belongs to the asymptotic confidence regions whenever $\theta_1$ and $\theta_2$ belong to those regions. However, they are not necessarily symmetric about $\hat{f}_X(x)$ or $\hat{m}(x)$. To implement the asymptotic confidence regions, we need to compute $\EL_{f_X}(\theta;x)$ and $\EL_m(\theta;x)$. The Lagrange multiplier technique leads that the unique maximizing weights $w_i$ are $1/(n(1+\lambda_{f_X}F_{f_X}(Z_i,\theta;x)))$ and $1/(n(1+\lambda_mF_m(Z_i,Y_i,\theta;x)))$ for $\EL_{f_X}(\theta;x)$ and $\EL_m(\theta;x)$, respectively, where $\lambda_{f_X}\in\mathbb{R}$ and $\lambda_m\in\mathbb{R}$ are the solutions of
\begin{align*}
\sum_{i=1}^n\frac{F_{f_X}(Z_i,\theta;x)}{1+\lambda_{f_X}F_{f_X}(Z_i,\theta;x)}=0\quad\text{and}\quad\sum_{i=1}^n\frac{F_m(Z_i,Y_i,\theta;x)}{1+\lambda_mF_m(Z_i,Y_i,\theta;x)}=0.
\end{align*}

\section{Finite sample performance}\label{simulation}

\setcounter{equation}{0}
\setcounter{subsection}{0}

\subsection{Simulation study}

In this section, we show the results of two simulation studies. We conducted regression analysis on $\mbS^2$ with measurement errors since this problem is practically important and it is our main interest. In the first simulation study, we checked the estimation performance of our regression estimator $\hat{m}$. Since there exists no other method designed for this problem, we compared $\hat{m}$ with a regression estimator designed for the error-free case. In particular, we took the naive regression estimator $\hat{m}^{\rm naive}$ defined as (\ref{naive regression estimator}) with $X_i$ in (\ref{naive regression estimator}) being replaced by $Z_i$, to see the effect of using $K_{T_n}$ instead of $K^*_{T_n}$. We recall that $K^*_{T_n}$ is introduced by \cite{Hendriks (1990)} for the error-free case. In the second simulation study, we compared the two types of asymptotic confidence intervals for $m$ that we constructed in Theorems \ref{asymptotic distribution regression toruses} and \ref{empirical likelihood regression}, namely the confidence interval based on the asymptotic normality (AN) and the confidence interval based on the empirical likelihood (EL), respectively. We recall that they are all currently available confidence intervals for this problem.

%with its density being given by $f_X(x)=\kappa\exp(\kappa\cdot\Theta\tran x)/(2\pi(\exp(\kappa)-\exp(-\kappa)))$

For both simulation studies, we generated $X$ from the von Mises-Fisher distribution on $\mbS^2$ with concentration parameter 0.1 and mean direction $(1,1,1)\tran/\sqrt{3}$. As for the distribution of $U$ on $SO(3)$, we took the Laplace distribution with $\lambda=0.5$ for the ordinary-smooth scenario (S1), the Gaussian distribution with $\lambda=0.5$ for the super-smooth scenario (S2) and the von Mises-Fisher distribution with $\lambda=2$ and $A=I_3$ for the log-super-smooth scenario (S3). The definitions of the Laplace, Gaussian and von Mises-Fisher distributions on  $SO(3)$ are given in Section \ref{error distribution}. We generated $Y$ from the model
\[
Y=(\cos\varphi_X+\sin\varphi_X)\sin\theta_X+\cos\theta_X+\epsilon,
\]
where $\varphi_X\in[0,2\pi)$ and $\theta_X\in[0,\pi)$ are the angles satisfying
\[
X=(\cos\varphi_X\sin\theta_X,\sin\varphi_X\sin\theta_X,\cos\theta_X)\tran,
\]
and $\epsilon$ is the normal random variable with mean zero and standard deviation 0.5. We chose $T_n$ based on a 5-fold cross-validation and repeatedly generated $\{(Y_i,U_iX_i):1\leq i\leq n\}$ with $n=250$ and $500$ for $R=200$ times.

In the first simulation study, we compared the integrated squared bias (ISB), integrated variance (IV) and integrated mean squared error (IMSE) defined by
\begin{align}\label{Criterion}
\begin{split}
\mbox{ISB}&=\int_{\mbS^2}\left(R^{-1}\sum_{r=1}^R\tilde{m}^{(r)}(x)-m(x)\right)^2d\nu(x),\\
\mbox{IV}&=R^{-1}\sum_{r=1}^R\int_{\mbS^2}\left(R^{-1}\sum_{s=1}^R\tilde{m}^{(s)}(x)-\tilde{m}^{(r)}(x)\right)^2d\nu(x),\\
\mbox{IMSE}&=\mbox{ISB+IV}=R^{-1}\sum_{r=1}^R\int_{\mbS^2}(\tilde{m}^{(r)}(x)-m(x))^2d\nu(x),
\end{split}
\end{align}
where $\tilde{m}^{(r)}(x)$ is either $\hat{m}(x)$ or $\hat{m}^{\rm naive}(x)$ obtained from the sample in the $r$th repeat for $1\leq r\leq R$. In the second simulation study, we computed the coverage rate $C_{1-\alpha}(x)$ and average length $L_{1-\alpha}(x)$ of $R$ confidence intervals of level $(1-\alpha)\times100\%$ for each $x\in\mathcal{G}$ and $\alpha\in\{0.05,0.1\}$, where $\mathcal{G}$ is a dense grid of $\mbS^2$. We then compared $|\mathcal{G}|^{-1}\sum_{x\in\mathcal{G}}C_{1-\alpha}(x)$ and $|\mathcal{G}|^{-1}\sum_{x\in\mathcal{G}}L_{1-\alpha}(x)$, where $|\mathcal{G}|$ denotes the cardinality of $\mathcal{G}$.

\begin{table*}[!t]
\caption{\footnotesize Integrated squared bias (ISB), integrated variance (IV) and integrated mean squared error (IMSE) of $\hat{m}$ for scenarios (S1)-(S3) based on $R=200$ Monte-Carlo samples.}
\label{sphere estimation}
{\footnotesize\centering
\begin{center}
\begin{tabular}{ c c c c c c c c }
\hline
&&\multicolumn{2}{c}{(S1)}&\multicolumn{2}{c}{(S2)}&\multicolumn{2}{c}{(S3)}\\
\hline
$n$ & Criterion & $\hat{m}$ & $\hat{m}^{\rm naive}$ & $\hat{m}$ & $\hat{m}^{\rm naive}$ & $\hat{m}$ & $\hat{m}^{\rm naive}$\\
\hline
       & ISB  & 0.04  & 1.40  & 0.04  & 0.61  & 0.07  & 0.95  \\ \cline{2-8}
 250   & IV   & 1.08  & 0.32  & 0.54  & 0.31  & 0.67  & 0.30  \\ \cline{2-8}
       & IMSE & 1.12  & 1.72  & 0.58  & 0.92  & 0.74  & 1.25  \\ \cline{1-8}
       & ISB  & 0.04  & 1.39  & 0.04  & 0.61  & 0.06  & 0.94  \\ \cline{2-8}
 500   & IV   & 0.39  & 0.15  & 0.27  & 0.16  & 0.34  & 0.17  \\ \cline{2-8}
       & IMSE & 0.43  & 1.54  & 0.31  & 0.77  & 0.39  & 1.11  \\ \hline
\end{tabular}
\end{center}
}
\end{table*}

\begin{table*}[!t]
\caption{\footnotesize Average coverage rate $|\mathcal{G}|^{-1}\sum_{x\in\mathcal{G}}C_{1-\alpha}(x)$ (Cov) and average length $|\mathcal{G}|^{-1}\sum_{x\in\mathcal{G}}L_{1-\alpha}(x)$ (Len) of $(1-\alpha)\times100\%$ confidence intervals of $m$ for scenario (S1) based on $R=200$ Monte-Carlo samples.}
\label{sphere interval}
{\small\centering
\begin{center}
\begin{tabular}{ c c c c c c }
\hline
&&\multicolumn{2}{c}{Cov}&\multicolumn{2}{c}{Len}\\
$(1-\alpha)$ & $n$ & AN & EL & AN & EL\\
\hline
0.9   & 250 & 0.87 & 0.86 & 0.77 & 0.84 \\ \cline{2-6}
      & 500 & 0.90 & 0.89 & 0.55 & 0.57 \\ \hline
0.95  & 250 & 0.91 & 0.91 & 0.92 & 1.01 \\ \cline{2-6}
      & 500 & 0.95 & 0.94 & 0.66 & 0.69 \\ \hline
\end{tabular}
\end{center}
}
\end{table*}

Table \ref{sphere estimation} shows the result of the first simulation study. The IMSE values demonstrate that $\hat{m}$ behaves better than $\hat{m}^{\rm naive}$. In particular, the ISB values of $\hat{m}$ are always much smaller than those of $\hat{m}^{\rm naive}$, which is explained by the unbiased scoring property of $\hat{m}$
as demonstrated in Proposition \ref{unbiased scoring}. While the errors of both estimators decrease as the sample size increases, the decreasing speed for $\hat{m}$ is much faster than that for $\hat{m}^{\rm naive}$. This suggests that our regression estimator is a reasonable estimator. Table \ref{sphere interval} shows the result of the second simulation study. It demonstrates that both methods generally produce higher coverage rates and narrower confidence intervals as the sample size increases. This suggests that both are reasonable methods. The table also reveals that the AN-based intervals have higher coverage rates and shorter lengths than the EL-based intervals. This indicates that the AN-based method can be a better option than the EL-based method. However, the latter is also a good alternative.

\subsection{Real data analysis}

We analyzed the dataset `sunspots\_births' in the R package `rotasym' (\cite{Garcia-Portugues et al. (2021)}). The dataset was analyzed in \cite{Garcia-Portugues et al. (2020)} to test the rotational symmetry of sunspots. Sunspots are temporary phenomena on the sun that appear as spots darker than the surrounding areas. Sunspot regions are cooler than the surrounding areas since the convection is blocked by the solar magnetic field flux. Sunspots are important sources in the study of solar activity and their number and positions affect the earth's long-term climate, telecommunications networks, aircraft navigation systems and spacecrafts, among others. Hence, it is important to study the distribution of sunspots.

Sunspots usually appear as a group. The dataset `sunspots\_births' contains $n=51,303$ groups of newly born sunspots measured in the years from 1872 to 2018. Each group observation contains the mean longitude and latitude of sunspots in that group. However, sunspots usually last only from a few hours to a few days, and they move across the surface of the sun. Also, the sizes of sunspots, known to have diameters ranging from 16km to 160,000km, keep changing during their lifespans. Due to these reasons combined with technical limitations of measuring devices, it is not easy to measure the exact birth locations of sunspots. Indeed, it is well known that sunspots area observation may contain measurement errors (e.g. \cite{Baranyi et al. (2001)}). Hence, we may assume that the observed mean longitudes and latitudes contain measurement errors. However, since the levels of the measurement errors could be not too high, we took the Laplace distribution on $SO(3)$ for the measurement error distribution and estimated the density of the birth locations of sunspots based on the deconvolution density estimator defined at (\ref{deconvolution density estimator}). We took the four values 0, $\sqrt{0.05}$, $\sqrt{0.1}$ and $\sqrt{0.15}$ for the distribution parameter $\lambda$, to see how the choice of the parameter affects the resulting density estimates. We note that the estimator with $\lambda=0$ corresponds to the naive density estimator that does not take into account measurement errors. We took the smoothing parameter $T_n$ minimizing the classical least squares cross-validation criterion (\cite{Rudemo (1982)}, \cite{Bowman (1984)}). This kind of comparison scheme was adopted in \cite{Efromovich (1997)} for contaminated $\mathbb{S}^1$-valued data. In this data analysis, we also included interval estimation studied in Theorems \ref{asymptotic distribution density toruses} and \ref{empirical likelihood density}. We note that asymptotic distributions and asymptotic confidence intervals for densities on $\mbS^2$ have not been studied in the literature of deconvolution density estimation on $\mbS^2$.

\begin{figure} %862 613 scale in R plot
\center
\includegraphics[width=0.45\textwidth]{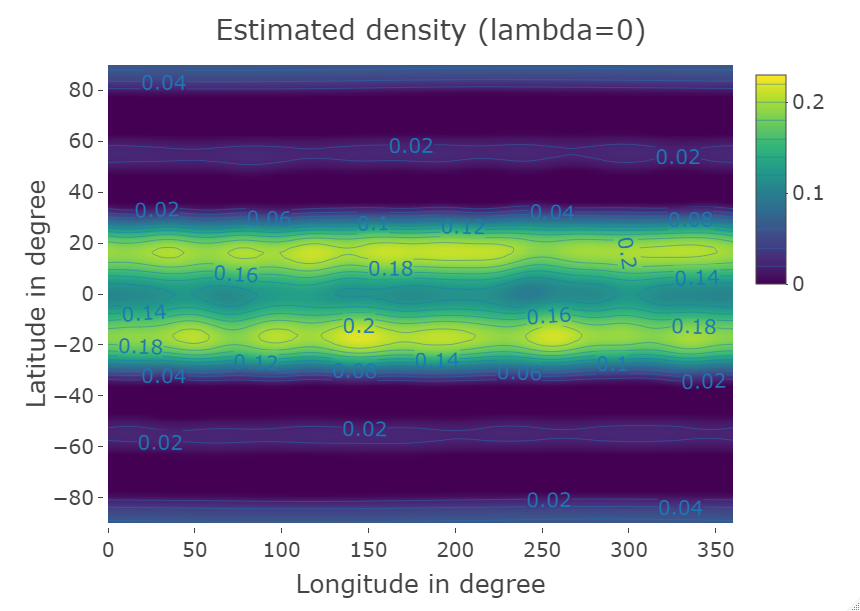}
\includegraphics[width=0.45\textwidth]{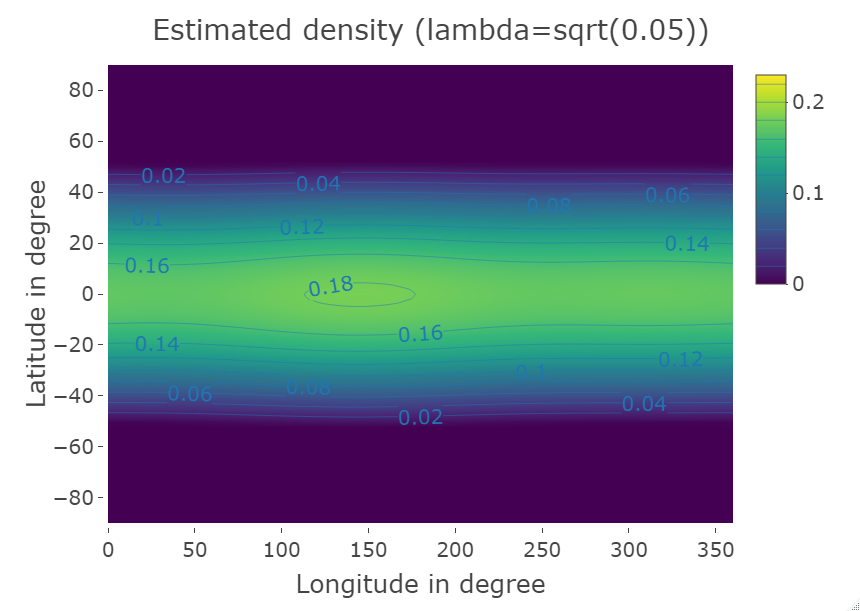}\\
\includegraphics[width=0.45\textwidth]{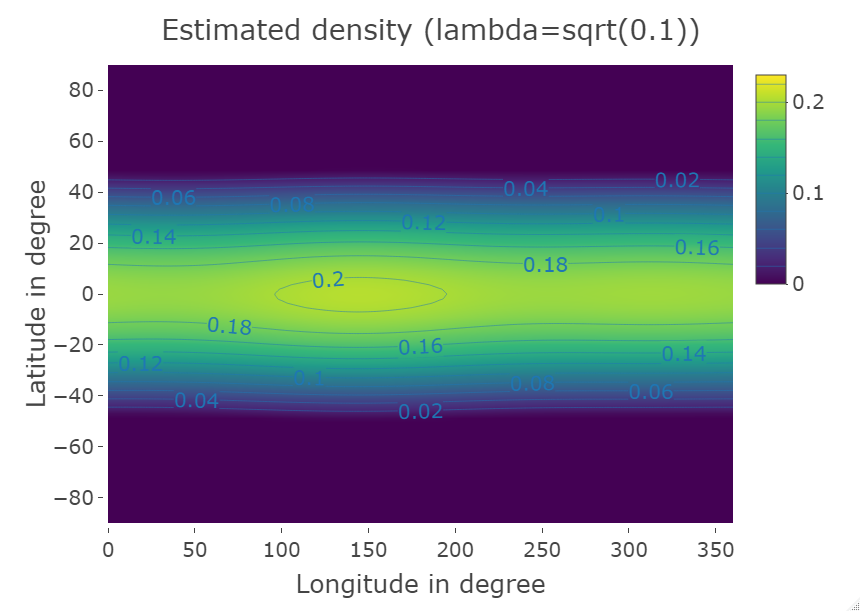}
\includegraphics[width=0.45\textwidth]{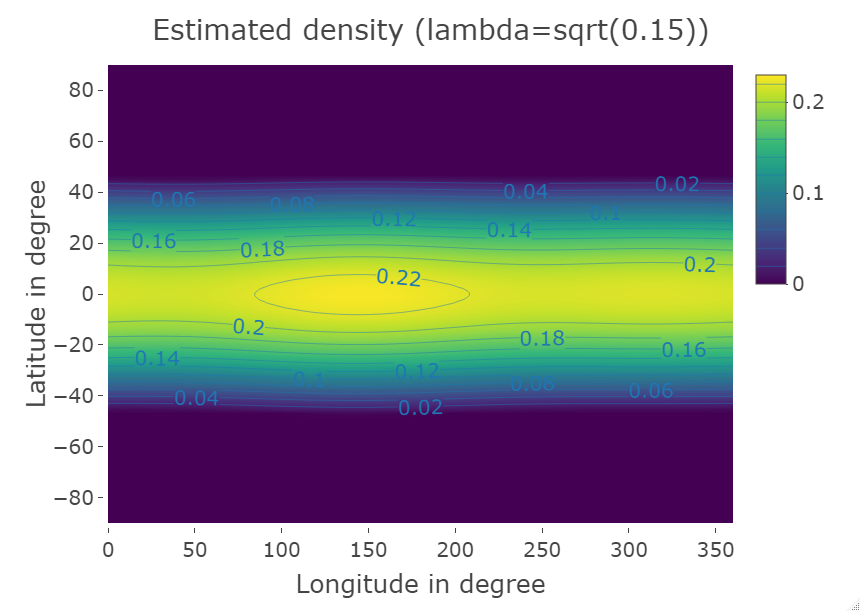}
\caption{\footnotesize The contour plots of the estimated densities based on the proposed method with $\lambda=0, \sqrt{0.05}, \sqrt{0.1}$ and $\sqrt{0.15}$. The color scale is the same for all plots.}
\label{estimated_density}
\end{figure}

\begin{figure}
\center
\includegraphics[width=0.45\textwidth]{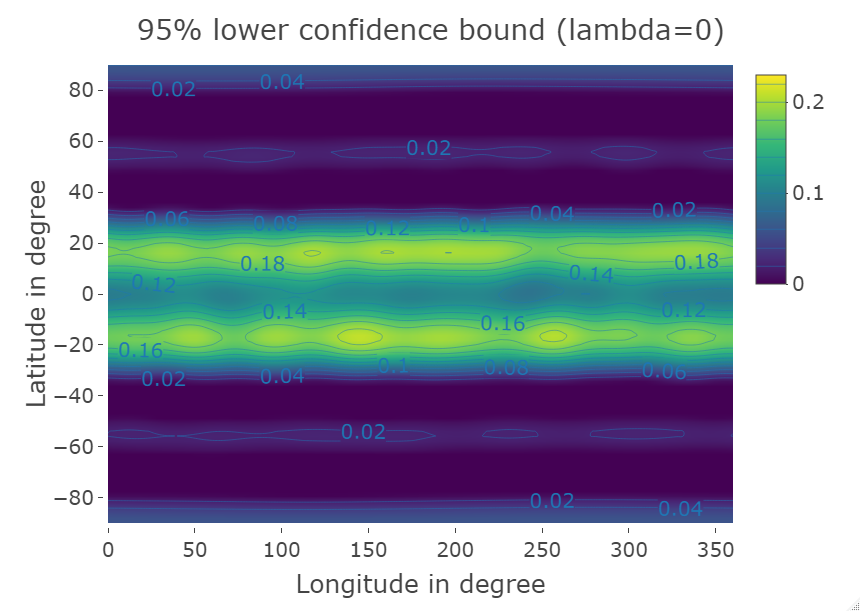}
\includegraphics[width=0.45\textwidth]{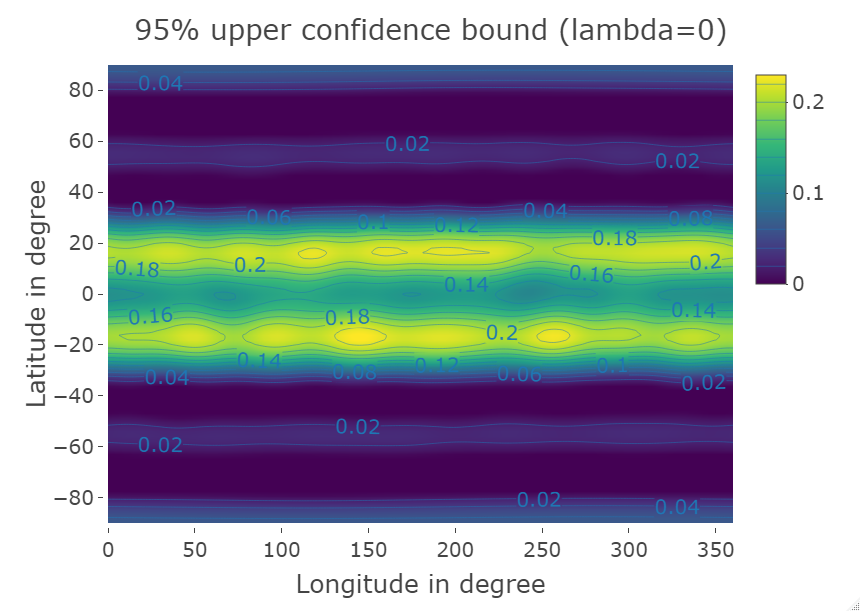}\\
\includegraphics[width=0.45\textwidth]{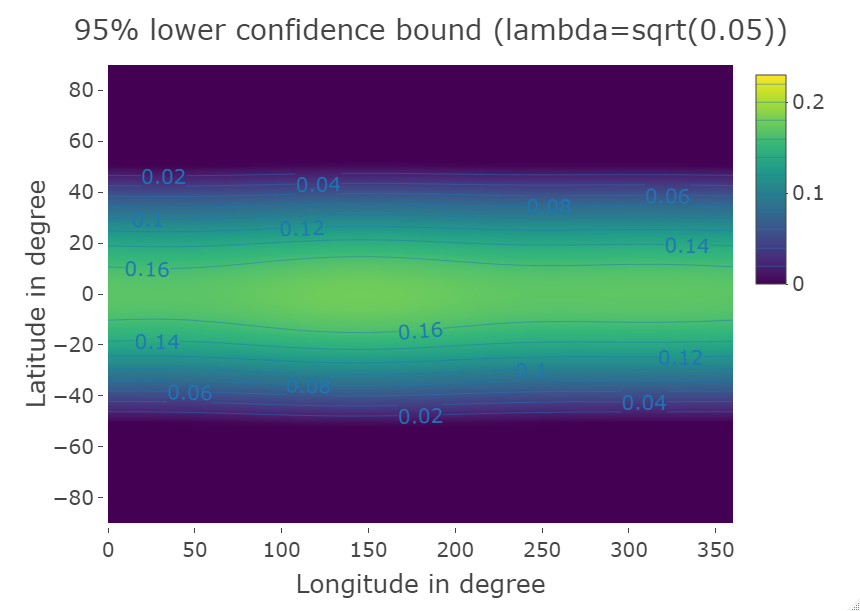}
\includegraphics[width=0.45\textwidth]{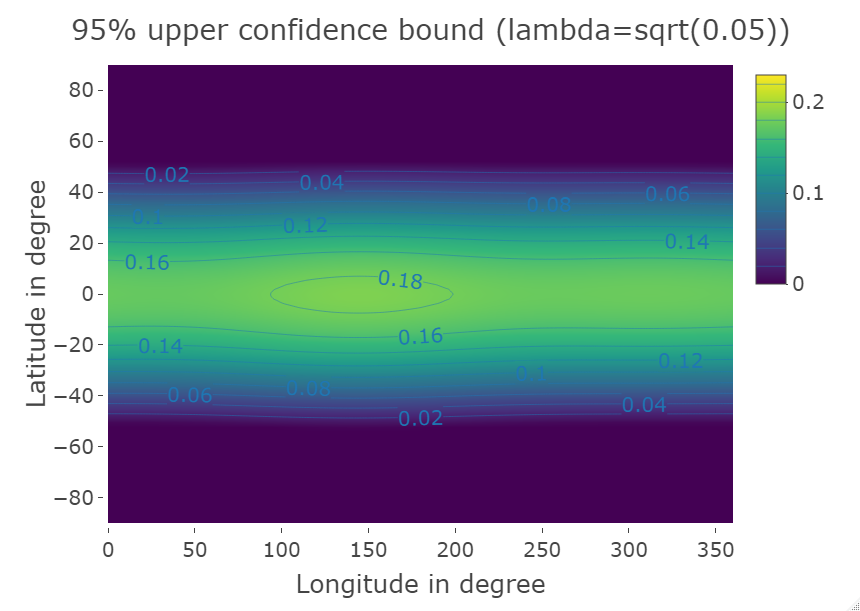}\\
\includegraphics[width=0.45\textwidth]{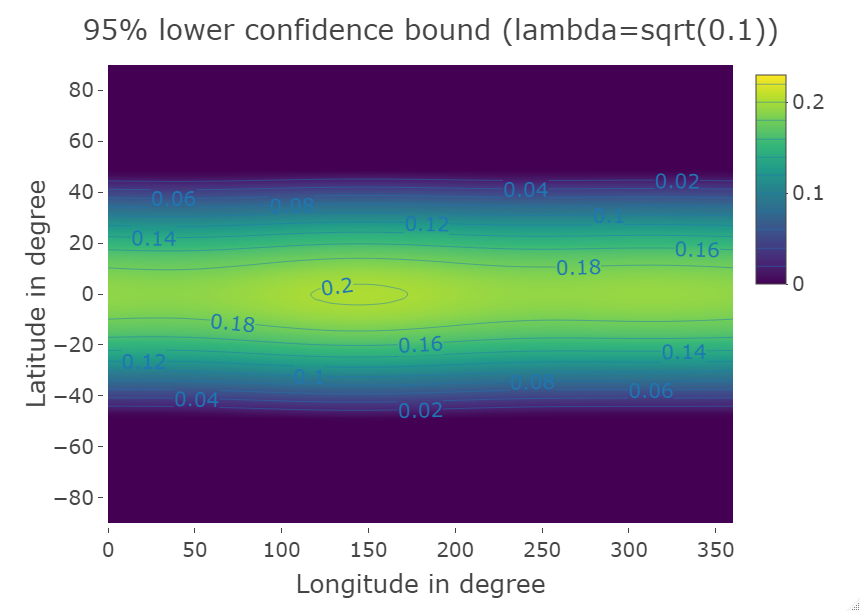}
\includegraphics[width=0.45\textwidth]{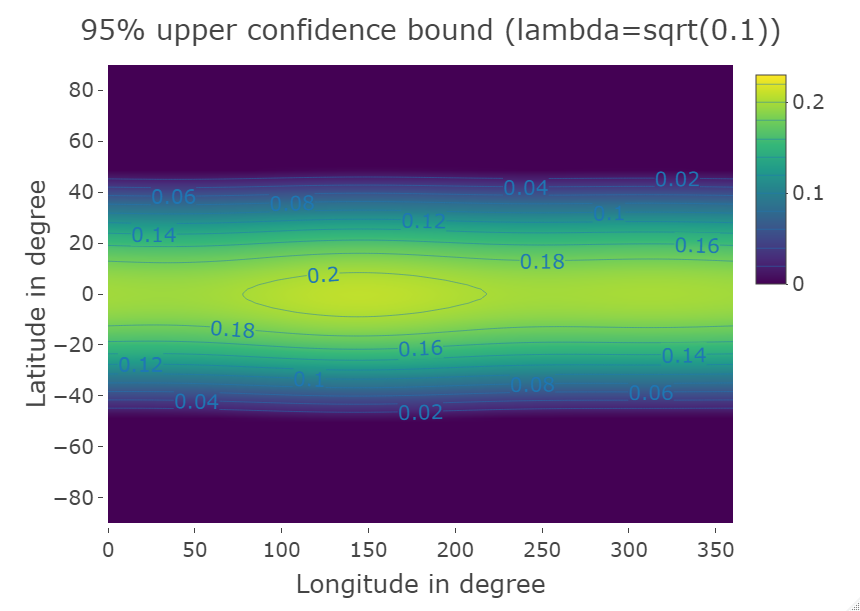}\\
\includegraphics[width=0.45\textwidth]{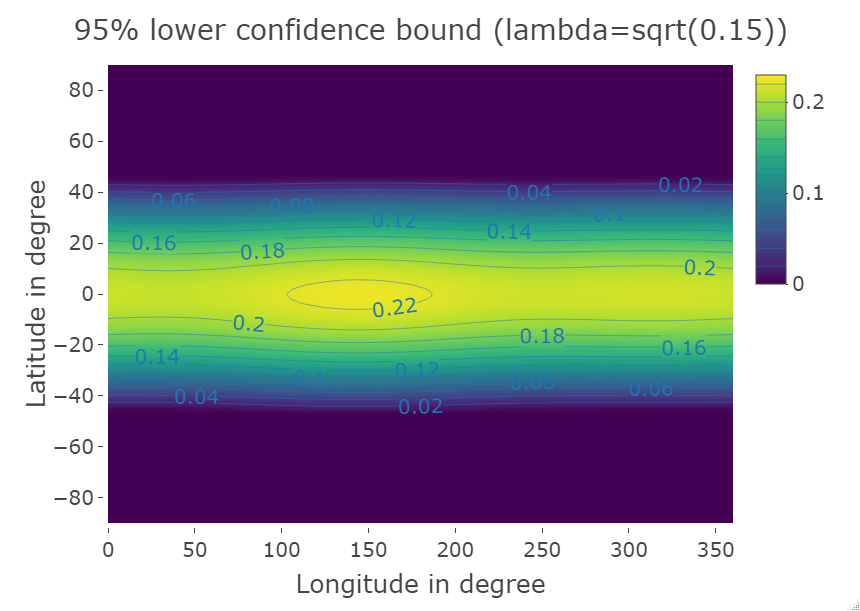}
\includegraphics[width=0.45\textwidth]{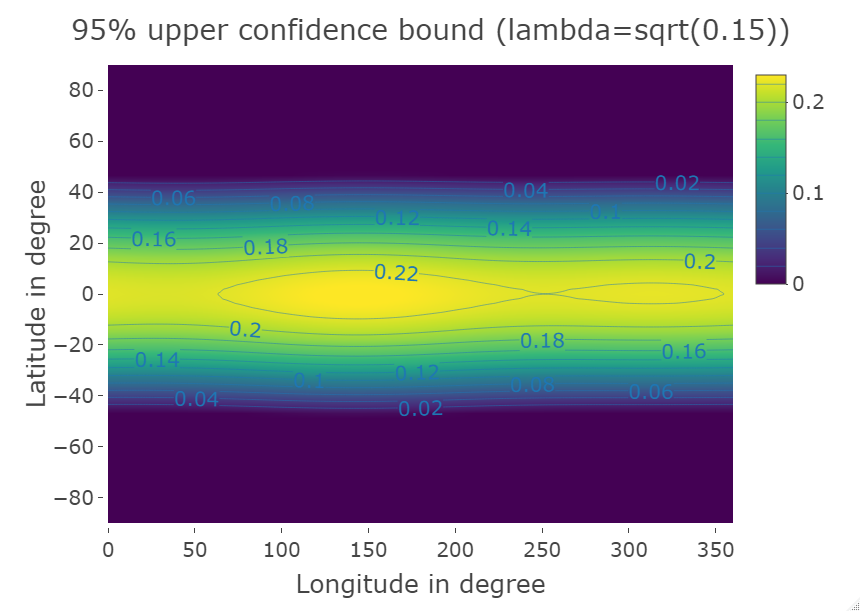}
\caption{\footnotesize The contour plots of the 95\% pointwise confidence intervals for the true density based on the asymptotic normality with $\lambda=0, \sqrt{0.05}, \sqrt{0.1}$ and $\sqrt{0.15}$. The color scale is the same for all plots.}
\label{AN CI_density}
\end{figure}

The contour plots of the estimated densities are depicted in Figure \ref{estimated_density}. The figure illustrates that, as $\lambda$ increases, the mass of the estimated density moves to the equator of the sun. It is well known that the rotating speed of the sun is the fastest at the equator and it decreases as the latitude goes up or down. Since sunspots are considered as a consequence of the twisted solar magnetic field caused by the fast rotating speed, the true density is likely to a higher mass as the latitude approaches to zero. It is also natural that the distribution of sunspots is symmetric about the equator and the density levels are horizontal due to the same reason. These justify the validity of the estimated densities for $\lambda>0$.

Figure \ref{AN CI_density} depicts the contour plots of the 95\% pointwise confidence intervals for the true density based on the asymptotic normality. The corresponding contour plots based on the empirical likelihood technique are omitted since they showed almost the same plots due to the large sample size. The upper confidence bounds on the right side of Figure \ref{AN CI_density} generally show wider peaks than the estimated densities in Figure \ref{estimated_density}, while the lower confidence bounds on the left side show opposite trends. Also, the length of each confidence interval is very short, which is very informative. We believe that these provide useful information in the analysis of sunspots.

\section*{Acknowledgements}
Research of Jeong Min Jeon was supported by the National Research Foundation of Korea (NRF) grant funded by the Korea government (MSIP) (No. 2020R1A6A3A03037314) and the European Research Council (2016-2021, Horizon 2020/ERC grant agreement No. 694409). Research of Ingrid Van Keilegom was supported by the European Research Council (2016-2021, Horizon 2020/ERC grant agreement No. 694409).

\newpage

\section*{}

\renewcommand{\theequation}{S.\arabic{equation}}
\renewcommand{\thesubsection}{S.\arabic{subsection}}
\renewcommand{\thelemma}{S.\arabic{lemma}}
\renewcommand{\theproposition}{S.\arabic{proposition}}
\setcounter{page}{1}

\setcounter{equation}{0}
\setcounter{subsection}{0}

\centerline{\bf \large Supplementary Material to}
\centerline{\bf \large `Density estimation and regression analysis on $\mathbb{S}^d$}
\centerline{\bf \large in the presence of measurement error'}
\centerline{\bf \large by Jeong Min Jeon and Ingrid Van Keilegom}

\bigskip

In the Supplementary Material, we provide some examples of the implementation of $D^l_{qr}(u)$ and $\tilde{\phi}^l_{qr}(f_U)$ for arbitrary $f_U$. We also provide all technical proofs. In the Supplementary Material, we denote by $ B^l(x)$ the $N(d,l)$-vector whose $q$th element equals $\Blq(x)$. We also let $\|\cdot\|_2$ denote the $L^2$-norm of $L^2((\mbS^d,\nu),\mathbb{C})$ and ${\rm (const.)}$ denote a generic positive constant.

\subsection{Implementation of $D^l_{qr}(u)$ and $\tilde{\phi}^l_{qr}(f_U)$ for arbitrary $f_U$}\label{phiU implementation}

\begin{itemize}
\item[1.] ($d=1$) It is well known that each $u\in SO(2)$ can be written as
\begin{align*}
        \begin{psmallmatrix}
        \cos\varphi_u & -\sin\varphi_u\\
        \sin\varphi_u & \cos\varphi_u
        \end{psmallmatrix}
\end{align*}
for some $\varphi_u\in[0,2\pi)$ and that $\int_{SO(2)}g(u)d\mu(u)=(2\pi)^{-1}\int_0^{2\pi}g(u)d\varphi_u$ for $g:SO(2)\rightarrow\mathbb{R}$. Using these and the definition of $\Blq$ given in Example \ref{spherical harmonics example}-1, we may show that
\begin{align*}
D^l_{11}(u)=\cos(l\varphi_u),&\quad D^l_{12}(u)=-\sin(l\varphi_u),\\
D^l_{21}(u)=\sin(l\varphi_u),&\quad D^l_{22}(u)=\cos(l\varphi_u).
\end{align*}
Using this, we have $\tilde{\phi}^l_{qr}(f_U)=(2\pi)^{-1}\int_0^{2\pi}f_U(u)D^l_{qr}(u)\,d\varphi_u$.
\item[2.] ($d=2$) We note that each $u\in SO(3)$ can be written as $R(\varphi_u)S(\theta_u)R(\psi_u)$ for some Euler angles $\varphi_u, \psi_u\in[0,2\pi)$ and $\theta_u\in[0,\pi)$, where
    \begin{align*}
    R(\vartheta)=
        \begin{psmallmatrix}
        \cos\vartheta & -\sin\vartheta & 0 \\
        \sin\vartheta & \cos\vartheta & 0\\
        0 & 0 & 1
        \end{psmallmatrix},
    \quad
    S(\vartheta)=
        \begin{psmallmatrix}
        \cos\vartheta & 0 & \sin\vartheta \\
        0 & 1 & 0\\
        -\sin\vartheta & 0 & \cos\vartheta
        \end{psmallmatrix}
    \end{align*}
for $\vartheta\in[0,2\pi)$ (Chapter 12.9 in Chirikjian (2012)). For $\Blq$ defined in Example \ref{spherical harmonics example}-2 and for $1\leq q,r\leq2l+1$, it holds that
\begin{align*}
D^l_{qr}(u)=e^{-\sqrt{-1}\cdot(q-l-1)\varphi_u}\cdot d^l_{qr}(\theta_u)\cdot e^{-\sqrt{-1}\cdot(r-l-1)\psi_u},
\end{align*}
where the definition of $d^l_{qr}(\theta_u)$ is given at (\ref{small d matrix}) (Chapter 12.9 in Chirikjian (2012)). Using this, we have
\begin{align*}
\tilde{\phi}^l_{qr}(f_U)=(8\pi^2)^{-1}\int_0^{2\pi}\int_0^\pi\int_0^{2\pi}f_U(u)D^l_{qr}(u)\sin\theta_u\,d\varphi_u\,d\theta_u\,d\psi_u;
\end{align*}
see Chapter 12.1 in Chirikjian (2012) for the representation of integration on $SO(3)$ with respect to the normalized Haar measure.
\end{itemize}

\subsection{Proof of Proposition \ref{deconvolution general}}

The proposition follows from
\begin{align*}
\phi^l_q(g\ast f)&=\int_{SO(d+1)}g(u)\int_{\mbS^d}f(u^{-1}x)\overline{\Blq(x)}\,d\nu(x)\,d\mu(u)\\
&=\int_{SO(d+1)}g(u)\int_{\mbS^d}f(x)\overline{\Blq(ux)}\,d\nu(x)\,d\mu(u)\\
&=\int_{SO(d+1)}g(u)\int_{\mbS^d}f(x)\sumr D^l_{qr}(u)\overline{\Blr(x)}\,d\nu(x)\,d\mu(u)\\
&=\sumr\int_{SO(d+1)}g(u)D^l_{qr}(u)\,d\mu(u)\int_{\mbS^d}f(x)\overline{\Blr(x)}\,d\nu(x)\\
&=\sumr\tilde{\phi}^l_{qr}(g)\phi^l_r(f),
\end{align*}
where we have used the rotation-invariant property of $\nu$ and that $\det(u)=1$ for the second equality, and (\ref{relation sum}) for the third equality.

\subsection{Proof of Proposition \ref{integral one}}

We first show that $\int_{\mbS^d}\Blq(x)d\nu(x)=0$ for $l>0$. Since $\{\Blq: l\in\mathbb{N}_0, 1\leq q\leq N(d,l)\}$ forms an orthonormal basis of $L^2((\mbS^d,\nu),\mathbb{C})$ and $ B^0_1\equiv(\nu(\mbS^d))^{-1/2}$, we have
\begin{align*}
\int_{\mbS^d}\Blq(x)d\nu(x)=(\nu(\mbS^d))^{1/2}\int_{\mbS^d}\Blq(x) B^0_1(x)d\nu(x)=0
\end{align*}
for $l>0$. Hence,
\begin{align*}
\int_{\mbS^d}K_{T_n}(x,z)d\nu(x)=&\suml\sumq\sumr(\tilde{\phi}^l(f_U))^{-1}_{qr}\,\overline{\Blr(z)}\int_{\mbS^d}\Blq(x)d\nu(x)\\
=&(\tilde{\phi}^0(f_U))^{-1}_{11}\,\overline{ B^0_1(z)}\int_{\mbS^d} B^0_1(x)d\nu(x)\\
=&1,
\end{align*}
where the last equality follows from the fact $\tilde{\phi}^0(f_U)=1$. This completes the proof.

\subsection{Proof of Proposition \ref{unbiased scoring}}

It suffices to show that
\begin{align*}
\E\left(\sumr(\tilde{\phi}^l(f_U))^{-1}_{qr}\overline{\Blr(Z)}\bigg|X\right)=\overline{\Blq(X)}.
\end{align*}
We note that
\begin{align*}
\overline{\Blr(Z)}=\overline{\Blr(UX)}=\sums D^l_{rs}(U)\overline{\Bls(X)},
\end{align*}
where the last equality follows from (\ref{deconvolution matrix form}). We also note that
\begin{align*}
\E\left(\sums D^l_{rs}(U)\overline{\Bls(X)}\bigg|X\right)=\sums\E(D^l_{rs}(U))\overline{\Bls(X)}=\sums\tilde{\phi}^l_{rs}(f_U)\overline{\Bls(X)},
\end{align*}
where the first equality follows from the assumption $U\perp X$. Hence, we have
\begin{align*}
\E\left(\sumr(\tilde{\phi}^l(f_U))^{-1}_{qr}\overline{\Blr(Z)}\bigg|X\right)&=\sums\overline{\Bls(X)}\sumr(\tilde{\phi}^l(f_U))^{-1}_{qr}\tilde{\phi}^l_{rs}(f_U)\\
&=\sums\overline{\Bls(X)}(I_{N(d,l)})_{qs}\\
&=\overline{\Blq(X)},
\end{align*}
which is the desired result.

\subsection{Proof of Proposition \ref{uniform consistency}}

We define $\tilde{f}_X(x)=n^{-1}\sum_{i=1}^nK_{T_n}(x,Z_i)$.
We note that $\hat{f}_X(x)=\RE(\tilde{f}_X(x))$ and
\begin{align*}
\E(\tilde{f}_X(x))&=\suml\sumq\Blq(x)\sumr(\tilde{\phi}^l(f_U))^{-1}_{qr}\phi^l_r(f_Z)\\
&=\suml\sumq\phi^l_q(f_X)\Blq(x),
\end{align*}
where the second equality follows from (\ref{relation sum}). Hence,
\begin{align*}
\sup_{x\in\mbS^d}|f_X(x)-\E(\tilde{f}_X(x))|=\sup_{x\in\mbS^d}\left|\sum_{l>T_n}\sumq\phi^l_q(f_X)\Blq(x)\right|=o(1),
\end{align*}
where the last equality follows from the assumption that the Fourier-Laplace series of $f_X$ converges uniformly. Also,
\begin{align*}
&\sup_{x\in\mbS^d}|\tilde{f}_X(x)-\E(\tilde{f}_X(x))|\\
\leq&\sup_{x\in\mbS^d}\suml\left|\sumq\Blq(x)\sumr(\tilde{\phi}^l(f_U))^{-1}_{qr}\left(n^{-1}\sum_{i=1}^n\overline{\Blr(Z_i)}-\phi^l_r(f_Z)\right)\right|\\
=&\sup_{x\in\mbS^d}\suml\left| B^l(x)\tran(\tilde{\phi}^l(f_U))^{-1}\left(n^{-1}\sum_{i=1}^n\overline{ B^l(Z_i)}-\phi^l(f_Z)\right)\right|\\
\leq&\sup_{x\in\mbS^d}\suml\| B^l(x)\|\cdot\left\|(\tilde{\phi}^l(f_U))^{-1}\left(n^{-1}\sum_{i=1}^n\overline{ B^l(Z_i)}-\phi^l(f_Z)\right)\right\|\\
\leq&\sup_{x\in\mbS^d}\suml\| B^l(x)\|\cdot\|(\tilde{\phi}^l(f_U))^{-1}\|_{\rm{op}}\left\|n^{-1}\sum_{i=1}^n\overline{ B^l(Z_i)}-\phi^l(f_Z)\right\|\\
=&\suml\sqrt{\frac{N(d,l)}{\nu(\mbS^d)}}\cdot\|(\tilde{\phi}^l(f_U))^{-1}\|_{\rm{op}}\left\|n^{-1}\sum_{i=1}^n\overline{ B^l(Z_i)}-\phi^l(f_Z)\right\|,
\end{align*}
where we have used the fact that $\| B^l(x)\|^2\equiv N(d,l)/\nu(\mbS^d)$ for the last equality. Hence,
\begin{align*}
&\E\left(\sup_{x\in\mbS^d}|\tilde{f}_X(x)-\E(\tilde{f}_X(x))|\right)\\
\leq&\suml\sqrt{\frac{N(d,l)}{\nu(\mbS^d)}}\cdot\|(\tilde{\phi}^l(f_U))^{-1}\|_{\rm{op}}\left(\sumr\E\left(\left|n^{-1}\sum_{i=1}^n\overline{\Blr(Z_i)}-\phi^l_r(f_Z)\right|^2\right)\right)^{1/2}\\
\leq&n^{-1/2}\suml\sqrt{\frac{N(d,l)}{\nu(\mbS^d)}}\cdot\|(\tilde{\phi}^l(f_U))^{-1}\|_{\rm{op}}\left(\sumr\E(|\Blr(Z)|^2)\right)^{1/2}\\
=&(\nu(\mbS^d))^{-1}n^{-1/2}\suml N(d,l)\|(\tilde{\phi}^l(f_U))^{-1}\|_{\rm{op}}.
\end{align*}

Now, we assume the case (S1)-(i). Then,
\begin{align*}
\suml N(d,l)\|(\tilde{\phi}^l(f_U))^{-1}\|_{\rm{op}}\leq1+c_1\sum_{l=1}^{[T_n]}N(d,l)l^\beta\leq{\rm (const.)}T_n^{\beta+d}
\end{align*}
since $N(d,l)=O(l^{d-1})$ as $l\rightarrow\infty$. By the choice (T1), it holds that
\begin{align*}
\sup_{x\in\mbS^d}|\tilde{f}_X(x)-f_X(x)|=\sup_{x\in\mbS^d}|\hat{f}_X(x)+\sqrt{-1}\cdot\IM(\tilde{f}_X(x))-f_X(x)|=o_p(1).
\end{align*}
Since $\sup_{x\in\mbS^d}|\hat{f}_X(x)-f_X(x)|\leq\sup_{x\in\mbS^d}|\tilde{f}_X(x)-f_X(x)|$,
the desired result follows. Now, we assume the case (S2)-(i). Then,
\begin{align*}
\suml N(d,l)\|(\tilde{\phi}^l(f_U))^{-1}\|_{\rm{op}}&\leq1+c_1\sum_{l=1}^{[T_n]}N(d,l)l^\alpha\exp(\gamma\cdot l^\beta)\\
&\leq{\rm (const.)}T_n^{\alpha+d}\exp(\gamma\cdot T_n^\beta).
\end{align*}
By the choice (T2), the result for the case (S2)-(i) similarly follows. Finally, we assume the case (S3)-(i). Then,
\begin{align*}
\suml N(d,l)\|(\tilde{\phi}^l(f_U))^{-1}\|_{\rm{op}}&\leq1+c_1\sum_{l=1}^{[T_n]}N(d,l)l^\alpha\exp(\gamma l^\beta(\log{l}-\xi_1))\\
&\leq{\rm (const.)}T_n^{\alpha+d}\exp(\gamma T_n^\beta(\log{T_n}-\xi_1)).
\end{align*}
Again by the choice (T3), the result for the case (S3)-(i) similarly follows. This completes the proof.

\subsection{Proof of Theorem \ref{L2 rates}}

We first prove the case of density estimation. Recall the definition of $\tilde{f}_X$ given in the proof of Proposition \ref{uniform consistency}. We note that
\begin{align*}
\E\left(\|\tilde{f}_X-f_X\|^2_2\right)=\E\left(\|\tilde{f}_X-\E(\tilde{f}_X)\|^2_2\right)+\|f_X-\E(\tilde{f}_X)\|^2_2.
\end{align*}
We first find the rate of $\|f_X-\E(\tilde{f}_X)\|^2_2$. We note that
\begin{align}\label{bias proof}
\begin{split}
&\|f_X-\E(\hat{f}_X)\|_2^2\\
&=\sum_{l>T_n}\sumq|\phi^l_q(f_X)|^2\\
&=\sum_{l>T_n}\sumq\lambda_l^{-2k}|(-1)^k\lambda_l^k\phi^l_q(f_X)|^2\\
&=(T_n(T_n+d-1))^{-2k}\sum_{l>T_n}\sumq(T_n(T_n+d-1))^{2k}\lambda_l^{-2k}|\phi^l_q(\Delta_{\mbS^d}^k(f_X))|^2\\
&\leq(T_n(T_n+d-1))^{-2k}\|\Delta^k_{\mbS^d}(f_X)\|^2_2,
\end{split}
\end{align}
where $-\lambda_l=-l(l+d-1)$ are the eigenvalues of the Laplace-Beltrami operator $\Delta_{\mbS^d}$ on $C^2(\mbS^d)$ and $\Delta_{\mbS^d}^k$ is the composition of $\Delta_{\mbS^d}$ for $k$-times. Since $\Delta^k_{\mbS^d}(f_X)$ is a continuous function on $\mbS^d$, we have $\|\Delta^k_{\mbS^d}(f_X)\|^2_2<\infty$. This gives
\begin{align*}
\|f_X-\E(\tilde{f}_X)\|^2_2=O(T_n^{-4k}).
\end{align*}
Now, we find the rate of $\E\left(\|\tilde{f}_X-\E(\tilde{f}_X)\|^2_2\right)=\int_{\mbS^d}\Var\left(\tilde{f}_X(x)\right)f_X(x)d\nu(x)$.
We note that
\begin{align*}
\Var\left(\tilde{f}_X(x)\right)=n^{-1}\Var\left(K_{T_n}(x,Z)\right)\leq n^{-1}\E(|K_{T_n}(x,Z)|^2).
\end{align*}
Since $f_X$ is bounded, it suffices to find the rate of $\int_{\mbS^d}\E(|K_{T_n}(x,Z)|^2)d\nu(x)$. It equals
\begin{align*}
\E\left(\int_{\mbS^d}|K_{T_n}(x,Z)|^2d\nu(x)\right)=&\suml\E\left(\sumq\left|\sumr(\tilde{\phi}^l(f_U))^{-1}_{qr}\overline{\Blr(Z)}\right|^2\right)\\
\leq&\suml\left\|(\tilde{\phi}^l(f_U))^{-1}\right\|_{\rm op}^2\E\left(\| B^l(Z)\|^2\right)\\
=&(\nu(\mbS^d))^{-1}\suml N(d,l)\left\|(\tilde{\phi}^l(f_U))^{-1}\right\|_{\rm op}^2,
\end{align*}
where the first equality follows from the orthonormality of $\{\Blq:1\leq q\leq N(d,l)\}$.

In the case of (a),
\begin{align*}
\suml N(d,l)\left\|(\tilde{\phi}^l(f_U))^{-1}\right\|_{\rm op}^2\leq1+c_1^2\sum_{l=1}^{[T_n]}N(d,l)l^{2\beta}\leq {\rm (const.)}T_n^{2\beta+d}.
\end{align*}
Therefore, we obtain
\begin{align*}
\E\left(\|\tilde{f}_X-f_X\|^2_2\right)=O(T_n^{-4k}+n^{-1}T_n^{2\beta+d}).
\end{align*}
This implies that
\begin{align*}
\|\hat{f}_X-f_X\|^2_2\leq\|\tilde{f}_X-f_X\|^2_2=O_p(T_n^{-4k}+n^{-1}T_n^{2\beta+d}).
\end{align*}
In the case of (b), we note that
\begin{align*}
\suml N(d,l)\left\|(\tilde{\phi}^l(f_U))^{-1}\right\|_{\rm op}^2\leq&1+c_1^2\sum_{l=1}^{[T_n]}N(d,l)l^{2\alpha}\exp(2\gamma\cdot T_n^{\beta})\\
\leq&{\rm (const.)}T_n^{2\alpha+d}\exp(2\gamma\cdot T_n^{\beta}).
\end{align*}
Therefore, we obtain
\begin{align*}
\E\left(\|\tilde{f}_X-f_X\|^2_2\right)=O(T_n^{-4k}+n^{-1}T_n^{2\alpha+d}\exp(2\gamma\cdot T_n^\beta)).
\end{align*}
This implies that
\begin{align*}
\|\hat{f}_X-f_X\|^2_2\leq\|\tilde{f}_X-f_X\|^2_2=O_p(T_n^{-4k}+n^{-1}T_n^{2\alpha+d}\exp(2\gamma\cdot T_n^\beta)).
\end{align*}
In the case of (c), similar arguments show that
\begin{align*}
\|\hat{f}_X-f_X\|^2_2&=O_p(T_n^{-4k}+n^{-1}T_n^{2\alpha+d}\exp(2\gamma\cdot T_n^\beta(\log{T_n}-\xi_1))).
\end{align*}
This completes the proof for the case of density estimation.

We now turn to the case of regression estimation. We write
\begin{align*}
\widetilde{m\cdot f_X}(x)&=n^{-1}\sum_{i=1}^n\left(\suml\sumq\Blq(x)\sumr(\tilde{\phi}^l(f_U))^{-1}_{qr}\overline{\Blr(Z_i)}\right)Y_i,\\
F(x)&=m(x)f_X(x)-\E(\widetilde{m\cdot f_X}(x)).
\end{align*}
We note that $\E\left(\|\widetilde{m\cdot f_X}-m\cdot f_X\|^2_2\right)=\E\left(\|\widetilde{m\cdot f_X}-\E(\widetilde{m\cdot f_X})\|^2_2\right)+\|F\|^2_2$.
We first approximate $\|F\|^2_2$. We note that
\begin{align}\label{expectation formula 1}
\begin{split}
\E\left(\sumr(\tilde{\phi}^l(f_U))^{-1}_{qr}\overline{\Blr(Z)}Y\right)&=\E\left(\E\left(\sumr(\tilde{\phi}^l(f_U))^{-1}_{qr}\overline{\Blr(Z)}Y\bigg|X\right)\right)\\
&=\E\left(\E\left(\sumr(\tilde{\phi}^l(f_U))^{-1}_{qr}\overline{\Blr(Z)}\bigg|X\right)m(X)\right)\\
&=\E(m(X)\overline{\Blq(X)})\\
&=\phi^l_q(m\cdot f_X),
\end{split}
\end{align}
where the second equality follows from the underlying assumption $U\perp (X,\epsilon)$, and the third equality follows from the proof of Proposition \ref{unbiased scoring}.
From (\ref{expectation formula 1}), we have
\begin{align*}
\E\left(\widetilde{m\cdot f_X}(x)\right)&=\suml\sumq\phi^l_q(m\cdot f_X)\Blq(x).
\end{align*}
Since this is a partial sum of the Fourier-Laplace series of $m\cdot f_X$ at $x$, and $m\cdot f_X$ is $2k$-times continuously differentiable, by arguing as (\ref{bias proof}), we get
\begin{align}\label{bias order}
\|F\|_2^2=O(T_n^{-4k}).
\end{align}
We now approximate
\begin{align*}
\E\left(\|\widetilde{m\cdot f_X}-\E(\widetilde{m\cdot f_X})\|^2_2\right)=\int_{\mbS^d}\Var\left(\widetilde{m\cdot f_X}(x)\right)f_X(x)d\nu(x).
\end{align*}
We note that
\begin{align*}
\Var\left(\widetilde{m\cdot f_X}(x)\right)=&n^{-1}\Var\left(K_{T_n}(x,Z)Y\right)\\
\leq&n^{-1}\E(|K_{T_n}(x,Z)Y|^2)\\
\leq&n^{-1}\E(\E(|K_{T_n}(x,Z)|^2|X)\E(Y^2|X))\\
\leq&{\rm (const.)}n^{-1}\E(|K_{T_n}(x,Z)|^2),
\end{align*}
where the second inequality follows from the underlying assumption $U\perp (X,\epsilon)$, and the last inequality follows from the boundedness of $\E(Y^2|X=\cdot)$. Since $f_X$ is bounded, it suffices to find the rate of $\int_{\mbS^d}\E(|K_{T_n}(x,Z)|^2)d\nu(x)$. The rate is obtained for each smoothness scenario in the proof of the first part of the theorem.

Hence, in the case of (a), we have
\begin{align*}
\E\left(\|\widetilde{m\cdot f_X}-m\cdot f_X\|^2_2\right)=O(T_n^{-4k}+n^{-1}T_n^{2\beta+d}).
\end{align*}
This implies that
\begin{align*}
\|\widehat{m\cdot f_X}-m\cdot f_X\|^2_2\leq\|\widetilde{m\cdot f_X}-m\cdot f_X\|^2_2=O_p(T_n^{-4k}+n^{-1}T_n^{2\beta+d}),
\end{align*}
where $\widehat{m\cdot f_X}:=\RE(\widetilde{m\cdot f_X})$.
We note that
\begin{align*}
\inf_{x\in \mbS^d}|\hat{f}_X(x)|&\geq \inf_{x\in \mbS^d}(f_X(x)-|\hat{f}_X(x)-f_X(x)|)\\
&\geq\inf_{x\in \mbS^d}f_X(x)-\sup_{x\in \mbS^d}|\hat{f}_X(x)-f_X(x)|.
\end{align*}
This with Proposition \ref{uniform consistency} and the assumption $\inf_{x\in \mbS^d}f_X(x)>0$ entails that there exists a constant $c>0$ such that $\inf_{x\in \mbS^d}|\hat{f}_X(x)|\geq c$ with probability tending to one. Since
\begin{align*}
&\int_{\mbS^d}|\hat{m}(x)-m(x)|^2d\nu(x)\\
&=\int_{\mbS^d}\left|\frac{\widehat{m\cdot f_X}(x)}{\hat{f}_X(x)}-\frac{m(x)f_X(x)}{f_X(x)}\right|^2d\nu(x)\\
&\leq2\int_{\mbS^d}\frac{|\widehat{m\cdot f_X}(x)-m(x)f_X(x)|^2}{|\hat{f}_X(x)|^2}+\frac{(m(x))^2|\hat{f}_X(x)-f_X(x)|^2}{|\hat{f}_X(x)|^2}d\nu(x)\\
&\leq{\rm (const.)}(\|\widehat{m\cdot f_X}-m\cdot f_X\|^2_2+\|\hat{f}_X-f_X\|^2_2)
\end{align*}
with probability tending to one, the result for the case (a) follows. The cases of (b) and (c) similarly follow as in the case of (a). This completes the proof.

\subsection{Proof of Theorem \ref{asymptotic distribution density}}

We note that
\begin{align*}
\hat{f}_X(x)-f_X(x)=n^{-1}\sum_{i=1}^n(\RE(K_{T_n}(x,Z_i))-f_X(x)).
\end{align*}
We write $W_{ni}(x)=n^{-1}(\RE(K_{T_n}(x,Z_i))-f_X(x))$ and show that
\begin{align}\label{asymptotic distribution density to show 1}
\frac{\sum_{i=1}^nW_{ni}(x)-\E(\sum_{i=1}^nW_{ni}(x))}{\sqrt{\Var(\sum_{i=1}^nW_{ni}(x))}}\overset{d}{\longrightarrow}N(0,1).
\end{align}
For this, we check that the Lyapunov condition
\begin{align}\label{W lyapunov density}
\frac{\E(|W_{n1}(x)-\E(W_{n1}(x))|^{2+\varsigma})}{n^{\varsigma/2}(\Var(W_{n1}(x)))^{1+\varsigma/2}}\rightarrow0
\end{align}
holds for some constant $\varsigma>0$. In particular, we choose any $\varsigma>0$ for the cases of (S2)-(i) and (S3)-(i). For the case of (S1)-(i), we choose $\varsigma>0$ satisfying $p<\varsigma/((2d-q)\varsigma+2(d-q))$. Such $\varsigma$ exists since $\varsigma/((2d-q)\varsigma+2(d-q))\rightarrow1/(2d-q)$ as $\varsigma\rightarrow\infty$. (\ref{W lyapunov density}) is equivalent to
\begin{align}\label{V lyapunov density}
\frac{\E(|V_n(x)-\E(V_n(x))|^{2+\varsigma})}{n^{\varsigma/2}(\Var(V_n(x)))^{1+\varsigma/2}}\rightarrow0,
\end{align}
where $V_n(x)=\RE(K_{T_n}(x,Z))-f_X(x)$. Since the nominator in (\ref{V lyapunov density}) is bounded by 
\begin{align*}
{\rm (const.)}\E(|\RE(K_{T_n}(x,Z))|^{2+\varsigma}), 
\end{align*}
it suffices to show that
\begin{align}\label{for empirical part 1}
\frac{\E(|\RE(K_{T_n}(x,Z))|^{2+\varsigma})}{n^{\varsigma/2}(\Var(V_n(x)))^{1+\varsigma/2}}\rightarrow0.
\end{align}
Also, since
\begin{align}\label{for empirical part 2}
\begin{split}
|\E(V_n(x))|&\leq\left|\sum_{l>T_n}\sumq\phi^l_q(f_X)\Blq(x)\right|=o(1),\\
\E((V_n(x))^2)&=\E((\RE(K_{T_n}(x,Z)))^2)-f^2_X(x)+o(1),\\
\E((\RE(K_{T_n}(x,Z)))^2)&\rightarrow\infty,
\end{split}
\end{align}
it suffices to show that
\begin{align}\label{asymptotic distribution density to show 2}
\frac{\E(|\RE(K_{T_n}(x,Z))|^{2+\varsigma})}{n^{\varsigma/2}(\E((\RE(K_{T_n}(x,Z)))^2))^{1+\varsigma/2}}\rightarrow0.
\end{align}
We note that
\begin{align}\label{two plus delta upper bound}
\begin{split}
\E(|K_{T_n}(x,Z)|^{2+\varsigma})&=\int_{\mbS^d}|K_{T_n}(x,z)|^{2+\varsigma}f_Z(z)d\nu(z)\\
&\leq\sup_{x\in\mbS^d}f_X(x)\int_{\mbS^d}|K_{T_n}(x,z)|^{2+\varsigma}d\nu(z),
\end{split}
\end{align}
where the inequality follows from the fact $\sup_{z\in\mbS^d}f_Z(z)\leq\sup_{x\in\mbS^d}f_X(x)$. We note that
\begin{align}\label{K max}
\begin{split}
|K_{T_n}(x,z)|^\varsigma&\leq\left(\suml\| B^l(x)\|\|(\tilde{\phi}^l(f_U))^{-1}\|_{\rm op}\| B^l(z)\|\right)^\varsigma\\
&=\left((\nu(\mbS^d))^{-1}\suml N(d,l)\|(\tilde{\phi}^l(f_U))^{-1}\|_{\rm op}\right)^\varsigma.
\end{split}
\end{align}
We also note that
\begin{align}\label{K integral z}
\begin{split}
\int_{\mbS^d}|K_{T_n}(x,z)|^2d\nu(z)&=\suml\sumr\left|\sumq B^l_q(x)(\tilde{\phi}^l(f_U))^{-1}_{qr}\right|^2\\
&\leq\suml\| B^l(x)\|^2\|((\tilde{\phi}^l(f_U))^{-1})\tran\|^2_{\rm op}\\
&=(\nu(\mbS^d))^{-1}\suml N(d,l)\|(\tilde{\phi}^l(f_U))^{-1}\|^2_{\rm op},
\end{split}
\end{align}
where the first equality follows from the orthonormality of $\{\Blq:1\leq q\leq N(d,l)\}$. Combining (\ref{two plus delta upper bound}), (\ref{K max}) and (\ref{K integral z}), we have
\begin{align}\label{upper bounds}
\begin{split}
&\E(|K_{T_n}(x,Z)|^{2+\varsigma})\\
&\leq
\begin{cases}
{\rm (const.)}T_n^{(2+\varsigma)\beta+(1+\varsigma)d}, & \text{if (S1)-(i) holds}\\
{\rm (const.)}T_n^{(2+\varsigma)\alpha+(1+\varsigma)d}\exp((2+\varsigma)\gamma\cdot T_n^\beta), & \text{if (S2)-(i) holds}\\
{\rm (const.)}T_n^{(2+\varsigma)\alpha+(1+\varsigma)d}\exp((2+\varsigma)\gamma\cdot T_n^\beta(\log{T_n}-\xi_1)), & \text{if (S3)-(i) holds}.
\end{cases}
\end{split}
\end{align}
Since $\E(|\RE(K_{T_n}(x,Z))|^{2+\varsigma})\leq\E(|K_{T_n}(x,Z)|^{2+\varsigma})$, $\E(|\RE(K_{T_n}(x,Z))|^{2+\varsigma})$ attains the same upper bounds given in (\ref{upper bounds}). Using (B1)-(B3) with $\eta$ sufficiently close to 1 in the cases of (B2) and (B3), we obtain (\ref{asymptotic distribution density to show 2}). Therefore, we have (\ref{asymptotic distribution density to show 1}). This completes the proof.

\subsection{Proof of Theorem \ref{asymptotic distribution regression} and some remark}

We note that
\begin{align*}
\hat{m}(x)-m(x)&=\frac{1}{\hat{f}_X(x)}\frac{1}{n}\sum_{i=1}^n\RE(K_{T_n}(x,Z_i))(Y_i-m(x))\\
&=\frac{1}{f_X(x)}\frac{1}{n}\sum_{i=1}^n\RE(K_{T_n}(x,Z_i))(Y_i-m(x))\left(1+\frac{f_X(x)-\hat{f}_X(x)}{\hat{f}_X(x)}\right)\\
&=\sum_{i=1}^nW_{ni}(x)+\sum_{i=1}^nW_{ni}(x)\cdot \frac{f_X(x)-\hat{f}_X(x)}{\hat{f}_X(x)},
\end{align*}
where $W_{ni}(x)=(f_X(x))^{-1}n^{-1}\RE(K_{T_n}(x,Z_i))(Y_i-m(x))$. Hence,
\begin{align*}
&\frac{\hat{m}(x)-m(x)-\E(\sum_{i=1}^nW_{ni}(x))}{\sqrt{\Var(\sum_{i=1}^nW_{ni}(x))}}\\
&=\frac{\sum_{i=1}^nW_{ni}(x)-\E(\sum_{i=1}^nW_{ni}(x))}{\sqrt{\Var(\sum_{i=1}^nW_{ni}(x))}}+\frac{\sum_{i=1}^nW_{ni}(x)}{\sqrt{\Var(\sum_{i=1}^nW_{ni}(x))}}\cdot \frac{f_X(x)-\hat{f}_X(x)}{\hat{f}_X(x)}.
\end{align*}
Thus, it suffices to prove that
\begin{align}\label{First assertion}
\frac{\sum_{i=1}^nW_{ni}(x)-\E(\sum_{i=1}^nW_{ni}(x))}{\sqrt{\Var(\sum_{i=1}^nW_{ni}(x))}}\overset{d}{\longrightarrow}N(0,1)
\end{align}
and
\begin{align}\label{Second assertion}
\frac{\sum_{i=1}^nW_{ni}(x)}{\sqrt{\Var(\sum_{i=1}^nW_{ni}(x))}}\cdot \frac{f_X(x)-\hat{f}_X(x)}{\hat{f}_X(x)}=o_p(1).
\end{align}

For (\ref{First assertion}), we check that the Lyapunov condition
\begin{align*}
\frac{\E(|W_{n1}(x)-\E(W_{n1}(x))|^{2+\delta})}{n^{\delta/2}(\Var(W_{n1}(x)))^{1+\delta/2}}\rightarrow0
\end{align*}
holds for $\delta$ in (B4). This is equivalent to verifying that
\begin{align}\label{V lyapunov}
\frac{\E(|V_n(x)-\E(V_n(x))|^{2+\delta})}{n^{\delta/2}(\Var(V_n(x)))^{1+\delta/2}}\rightarrow0,
\end{align}
where $V_n(x)=\RE(K_{T_n}(x,Z))(Y-m(x))$. Since the nominator in (\ref{V lyapunov}) is bounded by ${\rm (const.)}\E(|V_n(x)|^{2+\delta})$, it suffices to show that
\begin{align*}
\frac{\E(|V_n(x)|^{2+\delta})}{n^{\delta/2}(\Var(V_n(x)))^{1+\delta/2}}\rightarrow0.
\end{align*}
Also, since
\begin{align*}
\E(|V_n(x)|^{2+\delta})&\leq{\rm (const.)}\E(|\RE(K_{T_n}(x,Z))|^{2+\delta}),\\
\E((V_n(x))^2)&\geq{\rm (const.)}\E((\RE(K_{T_n}(x,Z)))^2)\rightarrow\infty,\\
\E(V_n(x))&=\E(\RE(K_{T_n}(x,Z))(m(X)-m(x)))=o(1),
\end{align*}
it suffices to show that
\begin{align*}
\frac{\E(|\RE(K_{T_n}(x,Z))|^{2+\delta})}{n^{\delta/2}(\E((\RE(K_{T_n}(x,Z)))^2))^{1+\delta/2}}\rightarrow0.
\end{align*}
But, this follows as in the proof of (\ref{asymptotic distribution density to show 2}).

For (\ref{Second assertion}), we note that
\begin{align}\label{for empirical part 3}
\begin{split}
&|f_X(x)-\E(\hat{f}_X(x))|\\
&\leq\sum_{l>T_n}\sumq|\phi^l_q(f_X)\Blq(x)|\\
&=\sum_{l>T_n}\sumq\lambda_l^{-k}|(-1)^k\lambda_l^k\phi^l(f_X)\Blq(x)|\\
&\leq(T_n(T_n-d+1))^{-k}\sum_{l>T_n}\sumq\lambda_l^{-k}(T_n(T_n-d+1))^k|\phi^l(\Delta_{\mbS^d}^k(f_X))\Blq(x)|\\
&\leq(T_n(T_n-d+1))^{-k}\sum_{l>T_n}\sumq|\phi^l(\Delta_{\mbS^d}^k(f_X))\Blq(x)|\\
&=o(T_n^{-2k}),
\end{split}
\end{align}
where the last equality follows from the absolute convergence of the Fourier-Laplace series of $\Delta_{\mbS^d}^k(f_X)$.
Also, it holds that
\begin{align*}
\E(\hat{f}_X(x))-\hat{f}_X(x)=O_p(n^{-1/2}\cdot T_n^{\beta+d})
\end{align*}
as in the proof of Proposition \ref{uniform consistency}. This with (\ref{for empirical part 3}) implies that
\begin{align}\label{Pointwise rate density}
f_X(x)-\hat{f}_X(x)=o(T_n^{-2k})+O_p(n^{-1/2}\cdot T_n^{\beta+d})=o_p(1).
\end{align}
Hence, it suffices to show that
\begin{align}\label{Third assertion}
\frac{\E(\sum_{i=1}^nW_{ni}(x))}{\sqrt{\Var(\sum_{i=1}^nW_{ni}(x))}}\cdot(f_X(x)-\hat{f}_X(x))=o_p(1)
\end{align}
by (\ref{First assertion}) and the fact that $(\hat{f}_X(x))^{-1}=O_p(1)$. We note that
\begin{align}\label{for empirical part 5}
\begin{split}
&f_X(x)\cdot\left|\E\left(\sum_{i=1}^nW_{ni}(x)\right)\right|\\
=&|\E(\RE(K_{T_n}(x,Z))(Y-m(x)))|\\
\leq&|m(x)|\cdot\sum_{l>T_n}\sumq|\phi^l_q(f_X)\Blq(x)|+\sum_{l>T_n}\sumq|\phi^l_q(m\cdot f_X)\Blq(x)|\\
\leq&(T_n(T_n-d+1))^{-k}\bigg(|m(x)|\cdot\sum_{l>T_n}\sumq|\phi^l(\Delta_{\mbS^d}^k(f_X))\Blq(x)|\\
&+\sum_{l>T_n}\sumq|\phi^l(\Delta_{\mbS^d}^k(m\cdot f_X))\Blq(x)|\bigg)\\
=&o(T_n^{-2k}),
\end{split}
\end{align}
where the last equality follows from the absolute convergence of the Fourier-Laplace series of $\Delta_{\mbS^d}^k(f_X)$ and of $\Delta_{\mbS^d}^k(m\cdot f_X)$.
We also note that
\begin{align}\label{only ordinary core}
\Var\left(\sum_{i=1}^nW_{ni}(x)\right)&\geq{\rm (const.)}n^{-1}\E((\RE(K_{T_n}(x,Z)))^2)\geq{\rm (const.)}n^{-1}T_n^{2\beta+q}.
\end{align}
Thus, (\ref{Third assertion}) follows if
\begin{align}\label{Remaining assertion}
n^{1/2}T_n^{-(\beta+q/2)}T_n^{-4k}=O(1)\quad\text{and}\quad T_n^{d-q/2}T_n^{-2k}=O(1).
\end{align}
The first one at (\ref{Remaining assertion}) follows by (T1$'''$). The second one at (\ref{Remaining assertion}) follows by (A1)-(i) and (A1)-(ii) with $k>(2d-q)/4$. This completes the proof.

\begin{remark}\label{only ordinary}
We give details on why we do not cover the super-smooth and log-super-smooth scenarios. Under (S2)-(i)+(B2), we obtain the rate $o(T_n^{-2k})+O_p(n^{-1/2}T_n^{\alpha+d}\exp(\gamma\cdot T_n^\beta))$ in the place of the right hand side of the first equality at (\ref{Pointwise rate density}) and the lower bound ${\rm (const.)}n^{-1}T_n^{2\alpha+q}\exp(2\gamma(\eta\cdot T_n)^{\beta})$ in the place of the right hand side of the last inequality at (\ref{only ordinary core}). Hence,
we need
\begin{align}\label{First requirement}
n^{-1/2}T_n^{\alpha+d}\exp(\gamma\cdot T_n^\beta)=o(1)
\end{align}
for the last equality at (\ref{Pointwise rate density}) and need
\begin{align}\label{Second requirement}
n^{1/2}T_n^{-(\alpha+q/2)}\exp(-\gamma(\eta\cdot T_n)^{\beta})\cdot T_n^{-2k}\cdot (T_n^{-2k}+n^{-1/2}T_n^{\alpha+d}\exp(\gamma\cdot T_n^\beta))=O(1)
\end{align}
for (\ref{Third assertion}). For (\ref{First requirement}), we need a log-type speed for $T_n$, while (\ref{Second requirement}) does not hold with any log-type speed for $T_n$. The log-super-smooth scenario has a similar problem.
%We also obtain the last lower bound at (\ref{only ordinary core}) by $n^{-1}T_n^{2\alpha+q}\exp(2\gamma(\eta\cdot T_n)^{\beta})$ for (B2) and $n^{-1}T_n^{2\alpha+q}\exp(2\gamma(\eta\cdot T_n)^{\beta}(\log T_n-\zeta))$ for (B3).
\end{remark}

\subsection{Proof of Lemma \ref{rate of denominator complex}}

Using the condition (C) and the fact $\inf_{z\in\mbS^d}f_Z(z)\geq\inf_{x\in\mbS^d}f_X(x)$, we obtain
\begin{align*}
\E(|K_{T_n}(x,Z)|^2)&\geq\int_{\mbS^d}|K_{T_n}(x,z)|^2d\nu(z)\cdot\inf_{x\in\mbS^d}f_X(x)\\
&=\suml\|((\tilde{\phi}^l(f_U))^{-1})\tran\Bl(x)\|^2\cdot\inf_{x\in\mbS^d}f_X(x)\\
&\geq{\rm (const.)}\suml\|\Bl(x)\|^2(\sigma_{\rm min}((\tilde{\phi}^l(f_U))^{-1}))^2\\
&\geq{\rm (const.)}\suml\|\Bl(x)\|^2\|(\tilde{\phi}^l(f_U))^{-1}\|^2_{\op}\\
&\geq{\rm (const.)}\suml N(d,l)\|(\tilde{\phi}^l(f_U))^{-1}\|^2_{\op}.
\end{align*}
Hence, for each $0<\eta<1$, we have
\begin{align*}
&\E(|K_{T_n}(x,Z)|^2)\\
&\geq
\begin{cases}
{\rm (const.)}\suml N(d,l)\cdot l^{2\beta}, & \text{if (S1)-(ii) holds}\\
{\rm (const.)}\suml N(d,l)\cdot l^{2\alpha}\cdot \exp(2\gamma\cdot l^\beta), & \text{if (S2)-(ii) holds}\\
{\rm (const.)}\suml N(d,l)\cdot l^{2\alpha}\cdot \exp(2\gamma\cdot l^\beta(\log{l}-\xi_2)), & \text{if (S3)-(ii) holds}
\end{cases}\\
&\geq
\begin{cases}
{\rm (const.)}(\eta\cdot T_n)^{2\beta}(\suml N(d,l)-\sum_{l=0}^{[\eta\cdot T_n]+1}N(d,l)), & \text{if (S1)-(ii) holds}\\
{\rm (const.)}(\eta\cdot T_n)^{2\alpha}(\suml N(d,l)-\sum_{l=0}^{[\eta\cdot T_n]+1}N(d,l)) & \\
\hspace*{1.175cm}\cdot\exp(2\gamma\cdot(\eta\cdot T_n)^\beta), & \text{if (S2)-(ii) holds}\\
{\rm (const.)}(\eta\cdot T_n)^{2\alpha}(\suml N(d,l)-\sum_{l=0}^{[\eta\cdot T_n]+1}N(d,l)) & \\
\hspace*{1.175cm}\cdot\exp(2\gamma\cdot(\eta\cdot T_n)^\beta(\log{(\eta\cdot T_n)}-\xi_2)), & \text{if (S3)-(ii) holds}.
\end{cases}
\end{align*}
We now prove that
\begin{align}\label{N sum order}
T_n^{-d}\left(\suml N(d,l)-\sum_{l=0}^{[\eta\cdot T_n]+1}N(d,l)\right)\rightarrow{\rm (const.)}.
\end{align}
We note that
\begin{align*}
N(d,l)=\left(2+{\frac{d-1}{l}}\right)\frac{1}{(d-1)!}\frac{\Gamma(l+d-1)}{\Gamma(l)}.
\end{align*}
From Tricomi and Erdelyi (1951), it is known that
\begin{align*}
\frac{\Gamma(l+d-1)}{\Gamma(l)}=l^{d-1}\left(1+\frac{(d-1)(d-2)}{2l}+O(l^{-2})\right).
\end{align*}
Hence, by Faulhaber's formula, we have
\begin{align*}
\suml\frac{\Gamma(l+d-1)}{\Gamma(l)}=\frac{1}{d}[T_n]^d+o([T_n]^d).
\end{align*}
This with simple algebra gives $T_n^{-d}\suml N(d,l)\rightarrow2/(d!)$ and hence (\ref{N sum order}) follows. Thus, we get
\begin{align}\label{lower bounds}
\begin{split}
&\E(|K_{T_n}(x,Z)|^2)\\
&\geq
\begin{cases}
{\rm (const.)}T_n^{2\beta+d}, & \text{if (S1)-(ii) holds}\\
{\rm (const.)}T_n^{2\alpha+d}\exp(2\gamma\cdot(\eta\cdot T_n)^\beta), & \text{if (S2)-(ii) holds}\\
{\rm (const.)}T_n^{2\alpha+d}\exp(2\gamma\cdot(\eta\cdot T_n)^\beta(\log{T_n}+\log{\eta}-\xi_2)), & \text{if (S3)-(ii) holds}.
\end{cases}
\end{split}
\end{align}
This completes the proof.

\subsection{Proof of Lemma \ref{rate of denominator toruses}}

We first consider the case (G1). In this case, $K_{T_n}(x,z)$ is real-valued for all $x, z\in\mbS^1$ since $\Blq$ and $\tilde{\phi}^l_{qr}$ are real-valued. Hence, the result follows.

Now, we consider the case (G2). We note that
\begin{align*}
K_{T_n}(x,z)&=\suml s_l^{-1}\sum_{q=1}^{2l+1}\Blq(x)\overline{\Blq(z)}\\
&=\suml s_l^{-1}\frac{2l+1}{4\pi}\sum_{q=1}^{2l+1}(\cos((q-l-1)\varphi_x)+\sqrt{-1}\sin((q-l-1)\varphi_x))\\
&\qquad\quad\cdot(\cos((q-l-1)\varphi_z)-\sqrt{-1}\sin((q-l-1)\varphi_z))d^l_{q(l+1)}(\theta_x)d^l_{q(l+1)}(\theta_z).
\end{align*}
Hence,
\begin{align*}
&\RE(K_{T_n}(x,z))=\suml s_l^{-1}\frac{2l+1}{4\pi}\sum_{q=1}^{2l+1}\cos((q-l-1)(\varphi_x-\varphi_z))d^l_{q(l+1)}(\theta_x)d^l_{q(l+1)}(\theta_z),\\
&\IM(K_{T_n}(x,z))=\suml s_l^{-1}\frac{2l+1}{4\pi}\sum_{q=1}^{2l+1}\sin((q-l-1)(\varphi_x-\varphi_z))d^l_{q(l+1)}(\theta_x)d^l_{q(l+1)}(\theta_z).
\end{align*}
Thus,
\begin{align*}
&\int_{\mbS^2}(\RE(K_{T_n}(x,z)))^2d\nu(z)\\
&=\suml\sum_{l'=0}^{[T_n]}s_l^{-1}s_{l'}^{-1}\frac{(2l+1)(2l'+1)}{16\pi^2}\sum_{q=1}^{2l+1}\sum_{q'=1}^{2l'+1}d^l_{q(l+1)}(\theta_x)d^{l'}_{q'(l'+1)}(\theta_x)\\
&~~~~~~~~~~~~~~\cdot\int_0^{2\pi}\cos((q-l-1)(\varphi_x-\varphi_z))\cos((q'-l'-1)(\varphi_x-\varphi_z))d\varphi_z\\
&~~~~~~~~~~~~~~\cdot\int_0^\pi d^l_{q(l+1)}(\theta_z)d^{l'}_{q'(l'+1)}(\theta_z)\sin(\theta_z)d\theta_z,\\
&\int_{\mbS^2}(\IM(K_{T_n}(x,z)))^2d\nu(z)\\
&=\suml\sum_{l'=0}^{[T_n]}s_l^{-1}s_{l'}^{-1}\frac{(2l+1)(2l'+1)}{16\pi^2}\sum_{q=1}^{2l+1}\sum_{q'=1}^{2l'+1}d^l_{q(l+1)}(\theta_x)d^{l'}_{q'(l'+1)}(\theta_x)\\
&~~~~~~~~~~~~~~\cdot\int_0^{2\pi}\sin((q-l-1)(\varphi_x-\varphi_z))\sin((q'-l'-1)(\varphi_x-\varphi_z))d\varphi_z\\
&~~~~~~~~~~~~~~\cdot\int_0^\pi d^l_{q(l+1)}(\theta_z)d^{l'}_{q'(l'+1)}(\theta_z)\sin(\theta_z)d\theta_z.
\end{align*}
Since
\begin{align*}
&\int_0^{2\pi}\cos((q-l-1)(\varphi_x-\varphi_z))\cos((q'-l'-1)(\varphi_x-\varphi_z))d\varphi_z\\
&=
\begin{cases}
2\pi, & \text{if $q-l-1=q'-l'-1=0$},\\
\pi, & \text{if $q-l-1=q'-l'-1\neq0$},\\
0, & \text{else},
\end{cases}\\
&\int_0^{2\pi}\sin((q-l-1)(\varphi_x-\varphi_z))\sin((q'-l'-1)(\varphi_x-\varphi_z))d\varphi_z\\
&=
\begin{cases}
\pi, & \text{if $q-l-1=q'-l'-1\neq0$},\\
0, & \text{else},
\end{cases}
\end{align*}
we have
\begin{align*}
&\int_{\mbS^2}(\RE(K_{T_n}(x,z)))^2d\nu(z)-\int_{\mbS^2}(\IM(K_{T_n}(x,z)))^2d\nu(z)\\
&=\suml\sum_{l'=0}^{[T_n]}s_l^{-1}s_{l'}^{-1}\frac{(2l+1)(2l'+1)}{8\pi}d^l_{(l+1)(l+1)}(\theta_x)d^{l'}_{(l'+1)(l'+1)}(\theta_x)\\
&~~~~~~~~~~~~~~\cdot\int_0^\pi d^l_{q(l+1)}(\theta_z)d^{l'}_{q'(l'+1)}(\theta_z)\sin(\theta_z)d\theta_z.
\end{align*}
By equation (12) in Pagaran et al. (2006), it holds that
\begin{align*}
\int_0^\pi d^l_{(l+1)(l+1)}(\theta_z)d^{l'}_{(l'+1)(l'+1)}(\theta_z)\sin(\theta_z)d\theta_z=\frac{2}{2l+1}I(l=l').
\end{align*}
Thus, we have
\begin{align*}
&\int_{\mbS^2}(\RE(K_{T_n}(x,z)))^2d\nu(z)-\int_{\mbS^2}(\IM(K_{T_n}(x,z)))^2d\nu(z)\\
&=\suml s_l^{-2}(d^l_{(l+1)(l+1)}(\theta_x))^2\frac{2l+1}{4\pi}\geq0.
\end{align*}
Combining this with the proof of Lemma \ref{rate of denominator complex} entails that $\int_{\mbS^2}(\RE(K_{T_n}(x,z)))^2d\nu(z)$ achieves the lower bounds given in (\ref{lower bounds}). Then, by the fact $\inf_{z\in\mbS^d}f_Z(z)\geq\inf_{x\in\mbS^d}f_X(x)$ and the condition (A2)-(i), we get the desired result.

\subsection{Proof of Theorem \ref{asymptotic distribution density toruses}}

We prove the two assertions
\begin{align}\label{zero bias C.I. density}
\sqrt{n}\cdot\frac{\E(\RE(K_{T_n}(x,Z)))-f_X(x)}{\sqrt{\Var(\RE(K_{T_n}(x,Z)))}}=o(1)
\end{align}
and
\begin{align}\label{variance equivalence C.I. density}
\frac{\hat{s}_1(x)}{\sqrt{\Var(\RE(K_{T_n}(x,Z)))}}\overset{p}{\rightarrow}1.
\end{align}
Then, we get the desired result by combining (\ref{zero bias C.I. density}), (\ref{variance equivalence C.I. density}) and Theorem \ref{asymptotic distribution density}.

For the assertion (\ref{zero bias C.I. density}), we note that
\begin{align*}
\Var(\RE(K_{T_n}(x,Z)))&=\E((\RE(K_{T_n}(x,Z)))^2)-(\E(\RE(K_{T_n}(x,Z))))^2\\
&=\E((\RE(K_{T_n}(x,Z)))^2)-f^2_X(x)+o(1),\\
\E((\RE(K_{T_n}(x,Z)))^2)&\rightarrow\infty.
\end{align*}
Hence, it suffices to show that
\begin{align}\label{for empirical part 2.5}
\sqrt{n}\cdot\frac{\E(\RE(K_{T_n}(x,Z)))-f_X(x)}{\sqrt{\E((\RE(K_{T_n}(x,Z)))^2)}}=o(1).
\end{align}
We note that
\begin{align*}
\E(\RE(K_{T_n}(x,Z)))-f_X(x)=o(T_n^{-2k}),
\end{align*}
by (\ref{for empirical part 3}). Hence, (\ref{for empirical part 2.5}) follows from (\ref{for empirical part 3}) and Lemma \ref{rate of denominator toruses} if $\sqrt{n}\cdot T_n^{-(2k+\beta+d/2)}=O(1)$. But, the latter holds with the choice (T1$'$). Thus, the assertion (\ref{zero bias C.I. density}) follows.

For the assertion (\ref{variance equivalence C.I. density}), we show that
\begin{align}\label{converges in probability C.I. density}
\begin{split}
\frac{1}{n}\sum_{i=1}^n\RE(K_{T_n}(x,Z_i))&\overset{p}{\rightarrow}\E(\RE(K_{T_n}(x,Z))),\\
\frac{1}{n}\sum_{i=1}^n(\RE(K_{T_n}(x,Z_i)))^2&\overset{p}{\rightarrow}\E((\RE(K_{T_n}(x,Z)))^2).
\end{split}
\end{align}
For the first one at (\ref{converges in probability C.I. density}), it suffices to show that
\begin{align*}
\Var\left(\frac{1}{n}\sum_{i=1}^n\RE(K_{T_n}(x,Z_i))\right)=o(1).
\end{align*}
This follows since
\begin{align*}
\E((\RE(K_{T_n}(x,Z)))^2)=O(T_n^{2\beta+d})=O(n^{(2\beta+d)/(4k+2\beta+d)})=o(n).
\end{align*}
For the second one at (\ref{converges in probability C.I. density}), we apply Corollary 2 in Chapter 10 of Chow and Teicher (1997). Then, it suffices to show that
\begin{align}\label{weak law triangle}
\E\left(\frac{(\RE(K_{T_n}(x,Z)))^2}{\E((\RE(K_{T_n}(x,Z)))^2)}\cdot I\left(\frac{(\RE(K_{T_n}(x,Z)))^2}{\E((\RE(K_{T_n}(x,Z)))^2)}\geq n\varepsilon\right)\right)\rightarrow0
\end{align}
holds for any $\varepsilon>0$. One can prove that (\ref{weak law triangle}) holds using (\ref{asymptotic distribution density to show 2}) and Lemma \ref{rate of denominator toruses}. Thus, the assertion (\ref{variance equivalence C.I. density}) follows. This completes the proof.

\subsection{Proof of Theorem \ref{asymptotic distribution regression toruses}}

We check the two claims
\begin{align}\label{zero bias C.I.}
\sqrt{n}\cdot\frac{\E(\RE(K_{T_n}(x,Z))(Y-m(x)))}{\sqrt{\Var(\RE(K_{T_n}(x,Z))(Y-m(x)))}}=o(1)
\end{align}
and
\begin{align}\label{variance equivalence C.I.}
\frac{\hat{s}_2(x)}{\sqrt{\Var(\RE(K_{T_n}(x,Z))(Y-m(x)))}}\overset{p}{\rightarrow}1.
\end{align}
Then, we get the desired result by combining (\ref{zero bias C.I.}), (\ref{variance equivalence C.I.}) and Theorem \ref{asymptotic distribution regression}.

The claim (\ref{zero bias C.I.}) follows if we prove that
\begin{align}\label{for empirical part 4}
\sqrt{n}\cdot\frac{\E(\RE(K_{T_n}(x,Z))(Y-m(x)))}{\sqrt{\E((\RE(K_{T_n}(x,Z)))^2)}}=o(1),
\end{align}
since
\begin{align*}
&\Var(\RE(K_{T_n}(x,Z))(Y-m(x)))\\
&=\E((\RE(K_{T_n}(x,Z))(Y-m(x)))^2)-(\E(\RE(K_{T_n}(x,Z))(Y-m(x))))^2\\
&\geq{\rm (const.)}\E((\RE(K_{T_n}(x,Z)))^2)+o(1).
\end{align*}
We note that (\ref{for empirical part 4}) follows from (\ref{for empirical part 5}) and Lemma \ref{rate of denominator toruses} provided that $\sqrt{n}\cdot T_n^{-(2k+\beta+d/2)}=O(1)$. But, the latter holds with the choice (T1$'$). Thus, the claim (\ref{zero bias C.I.}) follows.

The claim (\ref{variance equivalence C.I.}) follows if we show that
\begin{align}\label{converges in probability C.I.}
\begin{split}
\frac{1}{n}\sum_{i=1}^n\RE(K_{T_n}(x,Z_i))(Y_i-\hat{m}(x))&\overset{p}{\rightarrow}\E(\RE(K_{T_n}(x,Z))(Y-m(x))),\\
\frac{1}{n}\sum_{i=1}^n(\RE(K_{T_n}(x,Z_i))(Y_i-\hat{m}(x)))^2&\overset{p}{\rightarrow}\E((\RE(K_{T_n}(x,Z))(Y-m(x)))^2).
\end{split}
\end{align}
For the first one at (\ref{converges in probability C.I.}), we note that
\begin{align*}
&\frac{1}{n}\sum_{i=1}^n\RE(K_{T_n}(x,Z_i))(Y_i-\hat{m}(x))\\
=&\frac{1}{n}\sum_{i=1}^n\RE(K_{T_n}(x,Z_i))(Y_i-m(x))+(m(x)-\hat{m}(x))\hat{f}_X(x)\\
=&\frac{1}{n}\sum_{i=1}^n\RE(K_{T_n}(x,Z_i))(Y_i-m(x))+\left(m(x)-\frac{m(x)f_X(x)+o_p(1)}{f_X(x)+o_p(1)}\right)\hat{f}_X(x)\\
=&\frac{1}{n}\sum_{i=1}^n\RE(K_{T_n}(x,Z_i))(Y_i-m(x))+o_p(1),
\end{align*}
where the second equality follows similarly as in the proof of Proposition \ref{uniform consistency}.
Now, since
\begin{align*}
n^{-1}\E((\RE(K_{T_n}(x,Z))(Y-m(x)))^2)\leq{\rm (const.)}n^{-1}\E((\RE(K_{T_n}(x,Z)))^2)=o(1),
\end{align*}
the first one at (\ref{converges in probability C.I.}) follows. For the second one at (\ref{converges in probability C.I.}), we note that
\begin{align*}
&\frac{1}{n}\sum_{i=1}^n(\RE(K_{T_n}(x,Z_i))(Y_i-\hat{m}(x)))^2\\
&=\frac{1}{n}\sum_{i=1}^n(\RE(K_{T_n}(x,Z_i))(Y_i-m(x)))^2+(m(x)-\hat{m}(x))^2\frac{1}{n}\sum_{i=1}^n(\RE(K_{T_n}(x,Z_i)))^2\\
&\quad+(m(x)-\hat{m}(x))\frac{2}{n}\sum_{i=1}^n(\RE(K_{T_n}(x,Z_i)))^2(Y_i-m(x))\\
&=\frac{1}{n}\sum_{i=1}^n(\RE(K_{T_n}(x,Z_i))(Y_i-m(x)))^2+o_p(1).
\end{align*}
Hence, it suffices to show that
\begin{align*}
\frac{1}{n}\sum_{i=1}^n(\RE(K_{T_n}(x,Z_i))(Y_i-m(x)))^2\overset{p}{\rightarrow}\E((\RE(K_{T_n}(x,Z))(Y-m(x)))^2).
\end{align*}
For this, we apply Corollary 2 in Chapter 10 of Chow and Teicher (1997). Then, it suffices to show that
\begin{align}\label{weak law triangle 2}
\begin{split}
\E\bigg(\frac{(\RE(K_{T_n}(x,Z))(Y-m(x)))^2}{\E((\RE(K_{T_n}(x,Z))(Y-m(x)))^2)}\cdot I\bigg(\frac{(\RE(K_{T_n}(x,Z))(Y-m(x)))^2}{\E((\RE(K_{T_n}(x,Z))(Y-m(x)))^2)}\geq n\varepsilon\bigg)\bigg)\rightarrow0
\end{split}
\end{align}
holds for any $\varepsilon>0$. One can prove that (\ref{weak law triangle 2}) holds using (B4), (\ref{asymptotic distribution density to show 2}) and Lemma \ref{rate of denominator toruses}. Thus, the claim (\ref{variance equivalence C.I.}) follows. This completes the proof.

\subsection{Proof of Theorem \ref{empirical likelihood density}}

For the proof, we apply Theorem 2.1 in Hjort et al. (2009). For this, we verify the conditions (A0)-(A3) in Hjort et al. (2009). Note that
\begin{align*}
\EL_{f_X}(\theta;x)=\max\left\{\prod_{i=1}^n(nw_i):w_i>0, \sum_{i=1}^nw_i=1, \sum_{i=1}^nw_iF^*_{f_X}(Z_i,\theta;x)=0\right\},
\end{align*}
where
\begin{align*}
F^*_{f_X}(Z_i,\theta;x)=\frac{n^{-1/2}F_{f_X}(Z_i,\theta;x)}{\sqrt{\Var(\RE(K_{T_n}(x,Z)))}}.
\end{align*}
We note that (A0) in Hjort et al. (2009) immediately follows from the condition (E1). For (A1) in Hjort et al. (2009), it suffices to show that
\begin{align}\label{emp lik density 1}
\frac{n^{-1/2}\sum_{i=1}^n(\RE(K_{T_n}(x,Z_i))-f_X(x))}{\sqrt{\Var(\RE(K_{T_n}(x,Z)))}}\overset{d}{\longrightarrow}N(0,1).
\end{align}
From (\ref{asymptotic distribution density to show 1}), we have
\begin{align*}
\frac{n^{-1/2}\sum_{i=1}^n(\RE(K_{T_n}(x,Z_i))-f_X(x))-n^{1/2}\E(\RE(K_{T_n}(x,Z))-f_X(x))}{\sqrt{\Var(\RE(K_{T_n}(x,Z)))}}\overset{d}{\longrightarrow}N(0,1).
\end{align*}
We also have
\begin{align*}
\sqrt{n}\cdot\frac{\E(\RE(K_{T_n}(x,Z)))-f_X(x)}{\sqrt{\Var(\RE(K_{T_n}(x,Z)))}}=o(1)
\end{align*}
by (\ref{zero bias C.I. density}). Combining the two results gives (\ref{emp lik density 1}). For (A2) in Hjort et al. (2009), it suffices to show that
\begin{align}\label{emp lik density 2}
\frac{n^{-1}\sum_{i=1}^n(\RE(K_{T_n}(x,Z_i))-f_X(x))^2}{\Var(\RE(K_{T_n}(x,Z)))}\overset{p}{\rightarrow}1.
\end{align}
Since
\begin{align*}
&n^{-1}\sum_{i=1}^n(\RE(K_{T_n}(x,Z_i))-f_X(x))^2\\
&=n^{-1}\sum_{i=1}^n(\RE(K_{T_n}(x,Z_i)))^2-2n^{-1}f_X(x)\sum_{i=1}^n\RE(K_{T_n}(x,Z_i))+f^2_X(x),\\
&\Var(\RE(K_{T_n}(x,Z)))\\
&=\E((\RE(K_{T_n}(x,Z)))^2)-(\E(\RE(K_{T_n}(x,Z))))^2,
\end{align*}
it suffices to show that
\begin{align}\label{for A2 density}
\begin{split}
\E(\RE(K_{T_n}(x,Z)))\rightarrow f_X(x),\\
n^{-1}\sum_{i=1}^n\RE(K_{T_n}(x,Z_i))\overset{p}{\rightarrow}\E(\RE(K_{T_n}(x,Z))),\\
n^{-1}\sum_{i=1}^n(\RE(K_{T_n}(x,Z_i)))^2\overset{p}{\rightarrow}\E((\RE(K_{T_n}(x,Z)))^2).
\end{split}
\end{align}
The first assertion at (\ref{for A2 density}) follows from (\ref{for empirical part 2}), and the second and last assertions at (\ref{for A2 density}) follow from (\ref{converges in probability C.I. density}). Hence, (\ref{emp lik density 2}) holds. For (A3) in Hjort et al. (2009), it suffices to show that
\begin{align}\label{emp lik density 3}
\max_{1\leq i\leq n}|F^*_{f_X}(Z_i,f_X(x);x)|\overset{p}{\rightarrow}0.
\end{align}
We note that, for any $\varepsilon>0$ and $\varsigma>0$,
\begin{align*}
&P\left(n^{-1/2}\max_{1\leq i\leq n}|\RE(K_{T_n}(x,Z_i))-f_X(x)|>\varepsilon\sqrt{\Var(\RE(K_{T_n}(x,Z)))}\right)\\
&\leq{\rm (const.)}\frac{\E(|\RE(K_{T_n}(x,Z))-f_X(x)|^{2+\varsigma})}{n^{\varsigma/2}(\Var(\RE(K_{T_n}(x,Z))))^{1+\varsigma/2}}\\
&\leq{\rm (const.)}\frac{\E(|\RE(K_{T_n}(x,Z))|^{2+\varsigma})}{n^{\varsigma/2}(\Var(\RE(K_{T_n}(x,Z))))^{1+\varsigma/2}}\\
&\rightarrow0,
\end{align*}
where the limit follows similarly as in the proof of (\ref{for empirical part 1}). Thus, (\ref{emp lik density 3}) holds. Now, Theorem 2.1 in Hjort et al. (2009) gives the desired result.

\subsection{Proof of Theorem \ref{empirical likelihood regression}}

We apply Theorem 2.1 in Hjort et al. (2009) to prove the theorem. We note that
\begin{align*}
\EL_m(\theta;x)=\max\left\{\prod_{i=1}^n(nw_i):w_i>0, \sum_{i=1}^nw_i=1, \sum_{i=1}^nw_iF^*_m(Z_i,Y_i,\theta;x)=0\right\},
\end{align*}
where
\begin{align*}
F^*_m(Z_i,Y_i,\theta;x)=\frac{n^{-1/2}F_m(Z_i,Y_i,\theta;x)}{\sqrt{\Var(\RE(K_{T_n}(x,Z))(Y-m(x)))}}.
\end{align*}
Since the condition (A0) in Hjort et al. (2009) immediately follows from the condition (E2), it suffices to show that
\begin{align}\label{three assertions regression}
\begin{split}
\frac{n^{-1/2}\sum_{i=1}^n\RE(K_{T_n}(x,Z_i))(Y_i-m(x))}{\sqrt{\Var(\RE(K_{T_n}(x,Z))(Y-m(x)))}}&\overset{d}{\longrightarrow}N(0,1),\\
\frac{n^{-1}\sum_{i=1}^n(\RE(K_{T_n}(x,Z_i))(Y_i-m(x)))^2}{\Var(\RE(K_{T_n}(x,Z))(Y-m(x)))}&\overset{p}{\rightarrow}1,\\
\max_{1\leq i\leq n}|F^*_m(Z_i,Y_i,m(x);x)|&\overset{p}{\rightarrow}0
\end{split}
\end{align}
to verify the conditions (A1)-(A3) of Theorem 2.1 in Hjort et al. (2009).

For the first assertion at (\ref{three assertions regression}), we note that
\begin{align*}
\frac{n^{-1/2}\sum_{i=1}^n\RE(K_{T_n}(x,Z_i))(Y_i-m(x))-n^{1/2}\E(\RE(K_{T_n}(x,Z))(Y-m(x)))}{\sqrt{\Var(\RE(K_{T_n}(x,Z))(Y-m(x)))}}\overset{d}{\longrightarrow}N(0,1).
\end{align*}
This follows from the proof of Theorem \ref{asymptotic distribution regression}. Also, it holds that
\begin{align*}
\sqrt{n}\cdot\frac{\E(\RE(K_{T_n}(x,Z))(Y-m(x)))}{\sqrt{\Var(\RE(K_{T_n}(x,Z))(Y-m(x)))}}=o(1)
\end{align*}
by (\ref{zero bias C.I.}). Combining the two results gives the first assertion at (\ref{three assertions regression}). The second assertion at (\ref{three assertions regression}) follows from the facts
\begin{align*}
&\Var(\RE(K_{T_n}(x,Z))(Y-m(x)))=\E((\RE(K_{T_n}(x,Z))(Y-m(x)))^2)+o(1),\\
&\E((\RE(K_{T_n}(x,Z))(Y-m(x)))^2)\rightarrow\infty,\\
&n^{-1}\sum_{i=1}^n(\RE(K_{T_n}(x,Z_i))(Y_i-m(x)))^2\overset{p}{\rightarrow}\E((\RE(K_{T_n}(x,Z))(Y-m(x)))^2).
\end{align*}
For the third assertion at (\ref{three assertions regression}), we note that, for any $\varepsilon>0$ and $\delta>0$ in (B4),
\begin{align*}
&P\bigg(n^{-1/2}\max_{1\leq i\leq n}|\RE(K_{T_n}(x,Z_i))(Y_i-m(x))|>\varepsilon\sqrt{\Var(\RE(K_{T_n}(x,Z))(Y-m(x)))}\bigg)\\
&\leq{\rm (const.)}\frac{\E(|\RE(K_{T_n}(x,Z))(Y-m(x))|^{2+\delta})}{n^{\delta/2}(\Var(\RE(K_{T_n}(x,Z))(Y-m(x))))^{1+\delta/2}}\\
&\rightarrow0,
\end{align*}
where the limit follows similarly as in the proof of Theorem \ref{asymptotic distribution regression}. Now, Theorem 2.1 in Hjort et al. (2009) gives the desired result.

\bigskip
\large
\noindent
\textbf{References for Supplementary Material}

\bigskip

\small
\noindent
1. Chirikjian, G. S. (2012). {\it Stochastic Models, Information Theory, and Lie Groups, Volume 2}. Birkh\"{a}user Basel.\\
2. Chow, Y. S. and Teicher, H. (1997). {\it Probability Theory: Independence, Interchangeability, Martingales}. Springer-Verlag New York.\\
3. Hjort, N. L., McKeague, I. W. and Van Keilegom, I. (2009). Extending the scope of empirical likelihood. {\it Annals of Statistics}, \textbf{37}, 1079-1111.\\
4. Pagaran, J., Fritzsche, S. and Gaigalas, G. (2006). Maple procedures for the coupling of angular momenta. IX. Wigner D-functions and rotation matrices, {\it Computer Physics Communications}, \textbf{174}, 616-630.\\
5. Tricomi, F. G. and Erd\'{e}lyi, A. (1951). The asymptotic expansion of a ratio of gamma functions. {\it Pacific Journal of Mathematics}, \textbf{1}, 133-142.


\begin{thebibliography}{11}
\bibitem[Atkinson and Han(2012)]{Atkinson and Han (2012)} Atkinson, K. and Han, W. (2012). {\it Spherical Harmonics and Approximations on the Unit Sphere: An Introduction}. Springer-Verlag Berlin Heidelberg.
\bibitem[Baranyi et al.(2001)]{Baranyi et al. (2001)} Baranyi, T., Gyori, L., Ludm\'{a}ny, A. and Coffey, H. E. (2001). Comparison of sunspot area data bases. {\it Monthly Notices of the Royal Astronomical Society}, \textbf{323}, 223-230.
\bibitem[Belomestny and Goldenshluger(2021)]{Belomestny and Goldenshluger (2021)} Belomestny, D. and Goldenshluger, A. (2021). Density deconvolution under general assumptions on the distribution of measurement errors, {\it Annals of Statistics}, \textbf{49}, 615-649.
\bibitem[Bertrand et al.(2019)]{Bertrand et al. (2019)} Bertrand, A., Van Keilegom, I. and Legrand, C. (2019). Flexible parametric approach to classical measurement error variance estimation without auxiliary data. {\it Biometrics}, \textbf{75}, 297-307.
\bibitem[Boente et al.(2014)]{Boente et al. (2014)} Boente, G., Gonz\'{a}lez-Manteiga, W. and Rodr\'{\i}guez, D. (2009). Goodness-of-fit test for directional data. {\it Scandinavian Journal of Statistics}, \textbf{41}, 259-275.
\bibitem[Bowman (1984)]{Bowman (1984)} Bowman, A. W. (1984). An alternative method of cross-validation for the smoothing of density estimates. {\it Biometrika}, \textbf{72}, 353-360.
\bibitem[Chakraborty and Vemuri(2019)]{Chakraborty and Vemuri (2019)} Chakraborty, R. and Vemuri, B. C. (2019). Statistics on the Stiefel manifold: theory and applications.
{\it Annals of Statistics}, \textbf{47}, 415-438.
\bibitem[Chang(1989)]{Chang (1989)} Chang, T. (1989). Spherical regression with errors in variables. {\it Annals of Statistics}, \textbf{17}, 293-306.
\bibitem[Chen and Van Keilegom(2009)]{Chen and Van Keilegom (2009)} Chen, S. X. and Van Keilegom, I. (2009). A review on empirical likelihood methods for regression. {\it Test}, \textbf{18}, 415-447.
\bibitem[Chirikjian(2012)]{Chirikjian (2012)} Chirikjian, G. S. (2012). {\it Stochastic Models, Information Theory, and Lie Groups, Volume 2}. Birkh\"{a}user Basel.
%\bibitem[Chow and Teicher(1997)]{Chow and Teicher (1997)} Chow, Y. S. and Teicher, H. (1997). {\it Probability Theory: Independence, Interchangeability, Martingales}. Springer-Verlag New York.
\bibitem[Cuesta-Albertos et al.(2009)]{Cuesta-Albertos et al. (2009)} Cuesta-Albertos, J. A., Cuevas, A. and Fraiman, R. (2009). On projection-based tests for directional and compositional data. {\it Statistic and Computing}, \textbf{19}, 367-380.
\bibitem[Dattner et al.(2016)]{Dattner et al. (2016)} Dattner, I., Rei{\ss}, M. and Trabs, M. (2016). Adaptive quantile estimation in deconvolution with unknown error distribution. {\it Bernoulli}, \textbf{22}, 143-192.
\bibitem[Delaigle(2014)]{Delaigle (2014)} Delaigle, A. (2014). Nonparametric kernel methods with errors-in-variables: constructing estimators, computing them, and avoiding common mistakes. {\it Australian and New Zealand Journal of Statistics}, \textbf{56}, 105-124.
\bibitem[Delaigle et al.(2009)]{Delaigle et al. (2009)} Delaigle, A., Fan, J. and Carroll, R. J. (2009). A design-adaptive local polynomial estimator for the errors-in-variables problem. {\it Journal of the American Statistical Association}, \textbf{104}, 348-359.
\bibitem[Delaigle et al.(2015)]{Delaigle et al. (2015)} Delaigle, A., Hall, P. and Jamshidi, F. (2015). Confidence bands in non-parametric errors-in-variables regression. {\it Journal of the Royal Statistical Society: Series B (Statistical Methodology)}, \textbf{21}, 169-184.
\bibitem[Delaigle et al.(2008)]{Delaigle et al. (2008)} Delaigle, A., Hall, P. and Meister, A. (2008). On deconvolution with repeated measurements. {\it Annals of Statistics}, \textbf{36}, 665-685.
\bibitem[Efthimiou and Frye(2014)]{Efthimiou and Frye (2014)} Efthimiou, C. and Frye, C. (2014). {\it Spherical harmonics in p dimensions}. World Scientific Publishing Co. Pte. Ltd.
\bibitem[Efromovich(1997)]{Efromovich (1997)} Efromovich, S. (1997). Density estimation for the case of supersmooth measurement error. {\it Journal of the American Statistical Association}, \textbf{92}, 526-535.
\bibitem[Fan(1991a)]{Fan (1991a)} Fan, J. (1991a). On the optimal rates of convergence for nonparametric deconvolution problems. {\it Annals of Statistics}, \textbf{19}, 1257-1272.
\bibitem[Fan(1991b)]{Fan (1991b)} Fan, J. (1991b). Asymptotic normality for deconvolution kernel density estimators. {\it Sankhya}, \textbf{53}, 97-110.
\bibitem[Fan and Truong(1993)]{Fan and Truong (1993)} Fan, J. and Truong, Y. K. (1993). Nonparametric regression with errors in variables. {\it Annals of Statistics}, \textbf{21}, 1900-1925.
\bibitem[Gao et al.(2015)]{Gao et al. (2015)} Gao, F., Huang, X.-Y., Jacobs, N. A. and Wang, H. (2015). Assimilation of wind speed and direction
observations: results from real observation experiments. {\it Tellus A: Dynamic Meteorology and Oceanography}, \textbf{67}, 27132.
\bibitem[Garc\'{\i}a-Portugu\'{e}s et al.(2013)]{Garcia-Portugues et al. (2013)} Garc\'{\i}a-Portugu\'{e}s, E., Crujeiras, R. M. and Gonz\'{a}lez-Manteiga, W. (2013). Kernel density estimation for directional-linear data. {\it Journal of Multivariate Analysis}, \textbf{121}, 152-175.
\bibitem[Garc\'{\i}a-Portugu\'{e}s et al.(2020)]{Garcia-Portugues et al. (2020)} Garc\'{\i}a-Portugu\'{e}s, E., Paindaveine, D., and Verdebout, T. (2020). On optimal tests for rotational symmetry against new classes of hyperspherical distributions. {\it Journal of the American Statistical Association}, \textbf{115}, 1873-1887.
\bibitem[Garc\'{\i}a-Portugu\'{e}s et al.(2021)]{Garcia-Portugues et al. (2021)} Garc\'{\i}a-Portugu\'{e}s, E., Paindaveine, D., and Verdebout, T. (2021). rotasym: Tests for Rotational Symmetry on the Hypershpere. R package version 1.1.0.
\bibitem[Garc\'{\i}a-Portugu\'{e}s et al.(2016)]{Garcia-Portugues et al. (2016)} Garc\'{\i}a-Portugu\'{e}s, E., Van Keilegom, I., Crujeiras, R. M. and Gonz\'{a}lez-Manteiga, W. (2016). Testing parametric models in linear-directional regression. {\it Scandinavian Journal of Statistics}, \textbf{43}, 1178-1191.
\bibitem[Hall et al.(1987)]{Hall et al. (1987)} Hall, P., Watson, G. S. and Cabrera, J. (1987). Kernel density estimation with spherical data. {\it Biometrika}, \textbf{74}, 751-762.
\bibitem[Healy et al.(1998)]{Healy et al. (1998)} Healy, D. M., Hendriks, H. and Kim, P. T. (1998). Spherical deconvolution. {\it Journal of Multivariate Analysis}, \textbf{67}, 1-22.
\bibitem[Hendriks(1990)]{Hendriks (1990)} Hendriks, H. (1990). Nonparametric estimation of a probability density on a Riemannian manifold using Fourier expansions. {\it Annals of Statistics}, \textbf{18}, 832-849.
\bibitem[Hjort et al.(2009)]{Hjort et al. (2009)} Hjort, N. L., McKeague, I. W. and Van Keilegom, I. (2009). Extending the scope of empirical likelihood. {\it Annals of Statistics}, \textbf{37}, 1079-1111.
\bibitem[Huckemann et al.(2010)]{Huckemann et al. (2010)} Huckemann, S., Kim, P. T., Koo, J.-Y. and Munk, A. (2010). M\"{o}bius deconvolution on the hyperbolic plane with application to impedance density estimation. {\it Annals of Statistics}, \textbf{38}, 2465-2498.
\bibitem[Jeon et al.(2021)]{Jeon et al. (2021)} Jeon, J. M., Park, B. U. and Van Keilegom, I. (2021). Additive regression for non-Euclidean responses and predictors. {\it Annals of Statistics}, \textbf{49}, 2611-2641.
\bibitem[Jeon et al.(2022)]{Jeon et al. (2022)} Jeon, J. M., Park, B. U. and Van Keilegom, I. (2022). Nonparametric regression on Lie groups with measurement errors. {\it Annals of Statistics (under revision).}
\bibitem[Johannes(2009)]{Johannes (2009)} Johannes, J. (2009). Deconvolution with unknown measurement error distribution. {\it Annals of Statistics}, \textbf{37}, 2301-2323.
\bibitem[Johannes and Schwarz(2013)]{Johannes and Schwarz (2013)} Johannes, J. and Schwarz, M. (2013). Adaptive circular deconvolution by model selection under unknown error distribution. {\it Bernoulli}, \textbf{19}, 1576-1611.
\bibitem[Katznelson(2004)]{Katznelson (2004)} Katznelson, Y. (2004). {\it An introduction to harmonic analysis}. Cambridge University Press.
\bibitem[Kalf(1995)]{Kalf (1995)} Kalf, H. (1995). On the expansion of a function in terms of spherical harmonics in arbitrary dimensions. {\it Bulletin of the Belgian Mathematical Society}, \textbf{2}, 361-380.
%\bibitem[Kerkyacharian(2011)]{Kerkyacharian (2011)} Kerkyacharian, G., Pham Ngoc, T. M. and Picard, D. (2011). Localized spherical deconvolution. {\it Annals of Statistics}, \textbf{39}, 1042-1068.
\bibitem[Kim(1998)]{Kim (1998)} Kim, P. T. (1998). Deconvolution density estimation on SO(N). {\it Annals of Statistics}, \textbf{26}, 1083-1102.
\bibitem[Kim(2000)]{Kim (2000)} Kim, P. T. (2000). On the Characteristic Function of the Matrix von Mises-Fisher Distribution with Application to SO(N)-Deconvolution. {\it In: Gin\'{e}, E., Mason, D. M., Wellner, J. A. (eds) High Dimensional Probability II. Progress in Probability, Volume 47}, Birkh\"{a}user, Boston, MA.
\bibitem[Kim and Koo(2002)]{Kim and Koo (2002)} Kim, P. T. and Koo, J.-Y. (2002). Optimal spherical deconvolution. {\it Journal of Multivariate Analysis}, \textbf{80}, 21-42.
\bibitem[Kim et al.(2004)]{Kim et al. (2004)} Kim, P. T., Koo, J.-Y. and Park, H. J. (2004). Sharp minimaxity and spherical deconvolution for super-smooth error distributions. {\it Journal of Multivariate Analysis}, \textbf{90}, 384-392.
%\bibitem[Kim et al.(2016)]{Kim et al. (2016)} Kim, P. T., Koo, J.-Y. and Pham Ngoc, T. M. (2016). Supersmooth testing on the sphere over analytic classes. {\it Journal of Nonparametric Statistics}, \textbf{28}, 84-115.
\bibitem[Kim and Richards(2001)]{Kim and Richards (2001)} Kim, P. T. and Richards, D. St. P. (2001). Deconvolution density estimation on compact Lie groups. {\it Contemporary Mathematics}, \textbf{287}, 155-171.
%\bibitem[Lacour and Pham Ngoc(2014)]{Lacour and Pham Ngoc (2014)} Lacour, C. and Pham Ngoc, T. M. (2014). Goodness-of-fit test for noisy directional data. {\it Bernoulli}, \textbf{20}, 2131-2168.
\bibitem[Le\'{o}n et al.(2006)]{Leon et al. (2006)} Le\'{o}n, C. A., Mass\'{e}, J.-C. and Rivest, L.-P. (2006). A statistical model for random rotations. {\it Journal of Multivariate Analysis}, \textbf{97}, 412-430.
\bibitem[Luo et al.(2011)]{Luo et al. (2011)} Luo, Z. M., Kim, P. T., Kim, T. Y. and Koo, J.-Y. (2011). Deconvolution on the Euclidean motion group SE(3). {\it Inverse Problems}, \textbf{27}, 035014.
%\bibitem[Mardia and Jupp(1999)]{Mardia and Jupp (1999)} Mardia, K. V. and Jupp, P. E. (1999). {\it Directional Statistics}. John Wiley \& Sons, Inc.
\bibitem[Marron and Alonso(2014)]{Marron and Alonso (2014)} Marron, J. S. and Alonso, A. M. (2014). Overview of object oriented data analysis. {\it Biometical Journal}, \textbf{5}, 732-753.
\bibitem[Meister(2009)]{Meister (2009)} Meister, A. (2009). {\it Deconvolution Problems in Nonparametric Statistics}. Springer-Verlag Berlin Heidelberg.
\bibitem[Nadarajah and Zhang(2017)]{Nadarajah and Zhang (2017)} Nadarajah, S. J. and Zhang, Y. (2017). Wrapped: An R package for circular data. {\it PLoS ONE}, \textbf{12}, e0188512.
\bibitem[Owen(2001)]{Owen (2001)} Owen, A. (2001). Empirical Likelihood. Chapman and Hall/CRC, London.
%\bibitem[Pagaran et al.(2006)]{Pagaran et al. (2006)} Pagaran, J., Fritzsche, S. and Gaigalas, G. (2006). Maple procedures for the coupling of angular momenta. IX. Wigner D-functions and rotation matrices, {\it Computer Physics Communications}, \textbf{174}, 616-630.
\bibitem[Pewsey and Garc\'{\i}a-Portugu\'{e}s(2021)]{Pewsey and Garcia-Portugues (2021)} Pewsey, A. and Garc\'{\i}a-Portugu\'{e}s, E. (2021). Recent advances in directional statistics. {\it Test}, In print.
\bibitem[Qui et al.(2014)]{Qui et al. (2014)} Qui, Y., Nordman, D. J. and Vardeman, S. B. (2014). A wrapped trivariate normal distribution and Bayes inference for 3-D rotations. {\it Statistica Sinica}, \textbf{24}, 897-917.
\bibitem[Rivest(1989)]{Rivest (1989)} Rivest, L. P. (1989). Spherical regression for concentrated Fisher-von Mises distributions. {\it Annals of Statistics}, \textbf{17}, 307-317.
\bibitem[Rosenthal et al.(2014)]{Rosenthal et al. (2014)} Rosenthal, M., Wu, W. U., Klassen, E. and Srivastava, A. (2014). Spherical regression models using projective
linear transformations. {\it Journal of the American Statistical Association}, \textbf{109}, 1615-1624.
\bibitem[Rudemo(1982)]{Rudemo (1982)} Rudemo, M. (1982). Empirical choice of histograms and kernel density estimators. {\it Scandinavian Journal of Statistics}, \textbf{9}, 65-78.
\bibitem[Sakurai and Napolitano(2017)]{Sakurai and Napolitano (2017)} Sakurai, J. J. and Napolitano, J. (2017). {\it Modern Quantum Mechanics}. Cambridge University Press.
\bibitem[Scealy and Welsh(2011)]{Scealy and Welsh (2011)} Scealy, J. L. and Welsh, A. H. (2011). Regression for compositional data by using distributions defined on the hypersphere. {\it Journal of the Royal Statistical Society: Series B (Statistical Methodology)}, \textbf{73}, 351-375.
\bibitem[Sei et al.(2013)]{Sei et al. (2013)} Sei, T., Shibata, H., Takemura, A., Ohara, K. and Takayama, N. (2013). Properties and applications of Fisher distribution on the rotation group. {\it Journal of Multivariate Analysis}, \textbf{116}, 440-445.
\bibitem[Stefanski and Carroll(1990)]{Stefanski and Carroll (1990)} Stefanski, L. A. and Carroll, R. J. (1990). Deconvolving kernel density estimators. {\it Statistics}, \textbf{21}, 169-184.
\bibitem[Terras(2013)]{Terras (2013)} Terras, A. (2013). {\it Harmonic Analysis on Symmetric Spaces - Euclidean Space, the Sphere, and the Poincar\'{e} Upper Half-Plane}. Springer-Verlag New York.
%\bibitem[Tricomi and Erd\'{e}lyi(1951)]{Tricomi and Erdelyi (1951)} Tricomi, F. G. and Erd\'{e}lyi, A. (1951). The asymptotic expansion of a ratio of gamma functions. {\it Pacific Journal of Mathematics}, \textbf{1}, 133-142.
\end{thebibliography}
\end{document}